\documentclass[12pt]{iopart}

\usepackage{iopams}

\usepackage{epsfig}
\usepackage{graphicx}
\usepackage{xcolor}
\usepackage{makeidx}
\usepackage{multicol}
\usepackage{subeqn}
\usepackage{dsfont}
\usepackage{enumitem}

\def\d{\mathrm{d}}

\def\eps{{\varepsilon}}
\def\R{\mathbb{R}}

\def\I{\mathbb{I}}
\def\O#1{{\mathcal O}\left(#1\right)}
\def\sech{{\, \mathrm{sech} }}

\newtheorem{theorem}{Theorem}
\newtheorem{lemma}{Lemma}
\newtheorem{corollary}{Corollary}
\newtheorem{hypothesis}{Hypothesis}  
\newtheorem{remark}{Remark}
\newtheorem{definition}{Definition}

\def\proof{\noindent {\em Proof: \, }}
\def\proofend{\hfill $ \Box $}

\begin{document}

\title[Impulsive pseudo-manifolds]{Impulsive perturbations to differential equations:
stable/unstable pseudo-manifolds, heteroclinic connections, and flux}

\author{Sanjeeva Balasuriya}

\address{School of Mathematical Sciences, University of Adelaide, Adelaide SA 5005, Australia}
\ead{sanjeevabalasuriya@yahoo.com}
\vspace{10pt}
\begin{indented}
\item[]\today
\end{indented}

\begin{abstract}
State-dependent time-impulsive perturbations to a two-dimensional autonomous flow with stable and unstable manifolds are analysed by posing in terms of an integral equation which is valid in both forwards- and backwards-time.  The impulses destroy the smooth invariant manifolds, necessitating new definitions for stable and unstable {\em pseudo}-manifolds.  Their time-evolution is characterised by solving a Volterra integral equation of the second kind with discontinuous inhomogeniety.  A criteria for heteroclinic trajectory persistence in  this impulsive context is developed, as is a quantification of an instantaneous flux across broken heteroclinic manifolds. Several examples, including a kicked Duffing oscillator and an underwater explosion in the vicinity of an eddy, are used to illustrate the theory.
\end{abstract}

%
\vspace{2pc}

\vspace{10pt}
\submitto{\NL}
%
\maketitle
%
%

\section{Introduction}
\label{sec:intro}

It is well-known that stable and unstable manifolds are global flow organisers in autonomous flows arising
from ordinary differential equations \cite{guckenheimerholmes,wiggins}.  These are time-varying in nonautonomous flows, and their evolution in
relation to one another once again has important transport consequences \cite{peacockhaller}.  For example, this is well-understood
in two-dimensional time-periodic \cite{romkedar,wiggins,guckenheimerholmes} or time-aperiodic \cite{aperiodic} flows, 
and one might attempt to optimise transport across \cite{optimal,mixer,l2mixer}, or control the location of
\cite{saddlecontrol,manifoldcontrol,controlnd}, such structures in fluidic applications.

If a differential equation is subject to an {\em impulse}, the pleasing phase-space structure necessary for defining
stable and unstable manifolds gets destroyed.  Trajectories are no longer continuous in time, and hence smooth manifolds
cannot be defined.  On the other hand such impulses offer a natural method for modelling certain
types of phenomena, such as under-sea eruptions/earthquakes, a missile or other object falling into a body of water, or 
the tapping of a microfluidic device to incite mixing.  How would the modification to the fluid velocity as a result of such an 
impulse influence stable and unstable manifolds which were previously present?  What is the impact on fluid transport?

Thinking of impulses as simply resetting trajectory locations is well-established in the more applied literature.  This attitude
enables one to think directly in an autonomous phase space, but with trajectories jumping to new locations at the impulse
times.  An intuitively pleasing application, for example, is in controlling trajectories (in chaotic or other regimes); once a
trajectory starts exhibiting `bad' behaviour (such as getting influenced by an unstable manifold or chaotic attractor and getting pulled away), one
can think of resetting it to a previous `good' location.  This would be through the imposition of an impulse.  After the trajectory once
again approaches the same `bad' behaviour location, the impulse can be reapplied, and so on, resulting in a periodic trajectory
forced by periodic controlling impulses.  This and related ideas are available in the control and stabilisation
\cite[e.g.]{morrisgrizzle,grizzle,jianglu,osipov,cookekroll,hanlu,bressan} and neuroscience \cite{markram,catlla} literature.  These approaches, though useful in their particular context, 
do not capture the stable and unstable manifolds.

Let us be more concrete in describing the issues.  If $ x \in \Omega $, an $ n$-dimensional open connected set, 
the initial {\em intuition} might be to consider systems of the form
\begin{equation}
\fl \dot{x} = f(x) + \eps \sum_{i=1}^n  g_i(x,t) \delta(t-t_i)
\label{eq:diracde}
\end{equation}
where $ \delta $ is the Dirac delta  `function,' and $  \left\{ t_1, t_2, \cdots, t_n \right\} $ is an increasing
set of finite time values  at which the impulses  occur.  It is not assumed that the $ t_i $ are equally spaced; the system
(\ref{eq:diracde}) is nonautonomous.  The functions $ f $ and $ g_i $ are assumed smooth, and $ \left| \eps \right| $ is small.  Permitting the $ g_i $ to have $ x $-dependence means that the impact of the impulses is not uniform across 
$ \Omega $.
If the system (\ref{eq:diracde}) when $ \eps = 0 $ possesses a saddle fixed point with stable and unstable
manifolds, is it possible to characterise appropriate analogues of these when $ \eps \ne 0 $?
Trajectories starting at any initial
condition $ x(\beta) $ with $ \beta < t_1 $ would evolve continuously till $ x(t_1^-) $, but 
then must jump to $ x(t_1^+) $.  This jump is apparently quantified by $ \eps g_1 \left( x(t_1),t_1 \right) $, which immediately
leads to confusion since $ x(t_1) $ is not well-defined.  Is taking the left-hand limit the appropriate approach?  Or the right?  Or a combination?     Hence, (\ref{eq:diracde}) as it stands forms an ill-defined flow on $ \Omega $, a fact which has been
highlighted by several authors in the past \cite{griffithswalborn,coutinho,catlla}.
This problem arises because
the effect of the impulse is spatially-dependent (sometimes referred to as `state-dependent impulses' \cite{liuballinger,belleyvirgilio_duffing,belleyvirgilio_lienard}), as would be reasonable in applications such as underwater explosions.
This issue does {\em not} arise if the $ g_i $ are independent of $ x $ (as in state-independent kicks \cite{ott}), or if the $ x $-dependence
is such that there is no ambiguity in the jump (for example where a jump in one spatial variable depends on a different spatial
variable which does not encounter a jump \cite{linyoungkick,linyoungshear,wangyoung}, a specification which
explicitly uses only $ x(t_i^-) $ in its state-dependence \cite{dubeaukarrakchou}, or under other special conditions
\cite{wangoksasoglu}).

One resolution to this is to pose an autonomous differential equation which gets reset according to an {\em explicit} rule at specified discrete times; this is an
established method for addressing `impulsive differential equations' \cite{bainovsimeonov,barreiravalls,barreira,pan,hanlu,linyoungshear,linyoungkick,wangyoung,fennersiegmund,ott,wangoksasoglu}.
Usually, this rule is specified in one direction of time \cite{barreiravalls,barreira,pan,hanlu,linyoungshear,wangyoung,fennersiegmund,dubeaukarrakchou,liuballinger}, because of several reasons.  First, the function
specifying the resetting of trajectories need not be invertible in general, {\em unless} determined via a regularisation of
impulses such as in \cite{catlla}.  Second, if considering countable impulses occurring at $ t_1 < t_2 < t_3 < \cdots $ where
$ t_i \rightarrow \infty $, then
while it makes sense to flow time forward from time $ t < t_1 $ in, say, trying to understand a stable manifold \cite{fennersiegmund,pan,barreiravalls} or in establishing existence of solutions \cite{dubeaukarrakchou}, flowing backwards in time ``from infinity'' is troublesome.
Existing results from this perspective include proofs of existence of either the stable or the unstable
manifold (not both) by characterising the persistence of exponential decay estimates for the associated variational equation
\cite{barreiravalls,pan,fennersiegmund}, or proofs of chaotic dynamics or bifurcations
\cite{linyoungkick,linyoungshear,wangyoung,ott}.  The functional
analytic approach in these methods \cite{barreiravalls,barreira,pan,fennersiegmund} does not enable a method for actually locating and describing
the stable manifold.  In this article, explicitly characterising the time-variation of {\em both} the stable and the unstable manifold
will be pursued.  This will be possible by recasting the impulsive differential equation as an
integral equation according to a certain interpretation, befitting the ability of representing impulses in terms of distributions \cite{catlla}.   Furthermore, given the irregular time-dependence of the problem, formulating this on a nonautonomous (augmented) $ \Omega \times \R $ phase space, appropriately restricted, is a natural approach.

When viewed in the $ \Omega \times \R $ phase space, impulsive differential equations have 
a strong connection to autonomous vector fields which are discontinuous
\cite[e.g.]{battellifeckan,calamaifranca,kukucka,kunzekupper,duzhang,yagasakipiecewise}.  The reason is that in either situation, the $ \Omega \times
\R $ augmented phase space is partitioned by codimension-$1$  hyper-surfaces representing discontinuities, and the
evolution is governed by exactly how one matches trajectories crossing these discontinuity surfaces.  However, the temporal
discontinuities, i.e., impulses, addressed in this article are special in that time is a privileged independent variable in the augmented phase space, 
whose evolution is always given by $ \dot{t} = 1 $.  Several recent spatially discontinuous studies \cite{battellifeckan,calamaifranca,kukucka,duzhang,yagasakipiecewise} do have connections to this article in that they share the goal  of determining conditions on heteroclinic connections, 
while also being in the spirit of Melnikov theory \cite{melnikov,guckenheimerholmes,wiggins,siam_book}.  

In Section~\ref{sec:integral}, issues related to formalising (\ref{eq:diracde}) in terms of the standard impulsive differential
equations approach are discussed, and an integral equation formulation is proposed.
 Section~\ref{sec:pseudomanifolds} then defines the impulsive analogues of the stable and unstable
manifolds.  Clearly, these cannot exist as {\em manifolds} any longer, since impulses will destroy their smoothness.  This necessitates the definition of stable and unstable {\em pseudo-manifolds}, and their time-variation is obtained 
 in Theorems~\ref{theorem:unstabledirac} and \ref{theorem:stabledirac} by solving a Volterra integral equation of the second kind with discontinuous inhomogeneity over an unbounded domain.
It should be mentioned that the approach imposes neither time-periodicity nor volume-preservation.
Section~\ref{sec:heteroclinic} develops a condition for a persistent heteroclinic connection under impulses via Theorem~\ref{theorem:heteroclinic}, while 
Section~\ref{sec:flux} adapts the concept of an instantaneous flux \cite{aperiodic} to quantify the transport across a heteroclinic
manifold broken due to impulsive perturbations.  The above-mentioned theories of pseudo-manifolds, heteroclinic persistence,
and flux are respectively illustrated by examples in Sections~\ref{sec:parabolic}, \ref{sec:duffing}, \ref{sec:eddy} and
\ref{sec:expanding}
which follow each relevant section.  In particular, Section~\ref{sec:duffing} examines an impulsively kicked Duffing oscillator, 
characterising solutions which both forwards and backwards asymptote to the rest state, while Section~\ref{sec:eddy} addresses
the impact on water retention within an oceanic eddy due to a nearby explosion.

To the author's knowledge, this article is the first attempt to describe the time-variation of the locations of the analogues of stable and unstable manifolds in flows subject to spatially-dependent time impulses.  The approach is geometric in nature, appealing to physical intuition in 
the augmented phase space and---in this first attempt---is restricted to $ \Omega $ being two-dimensional.  Thus, for example,
the stable pseudo-manifold would be a time-varying curve in $ \Omega $ which when advected in forwards time 
collapses to an unstable  fixed point.  Though describing these entities with geometric intuition, the
development is not merely formal, and takes into account rigorous distributional derivatives while ensuring that errors are 
higher-order.

\section{Integral equation formulation}
\label{sec:integral}

Consider the `conceptual equation' (\ref{eq:diracde}) for which a well-defined formulation is sought.   First, some conditions
on the functions will be stated.

\begin{hypothesis}[Unperturbed flow conditions]
\label{hyp:fdirac}
The unperturbed ($\eps = 0$) system (\ref{eq:diracde}) is associated with the conditions
\begin{itemize}
\item[(a)] $ f \in {\mathrm{C}}^2 \left( \Omega \right) $ with $ D f $ bounded in 
$ \Omega $, an open connected two-dimensional set;
\item[(b)] There exists $ a \in \Omega $ such that $ f(a) = 0 $ and $ D f(a) $ possesses a positive and a negative eigenvalue.
\end{itemize}
\end{hypothesis}

The implication of Hypothesis~\ref{hyp:fdirac} is that when $ \eps = 0 $, (\ref{eq:diracde}) may as well be considered
as a differential equation
\begin{equation}
\fl \dot{x} = f(x)
\label{eq:unpde}
\end{equation}
for $ x \in \Omega $,
in which $ a $ is a fixed point which has one-dimensional stable and unstable manifolds emanating from it.  This differential
formulation is not possible when $ \eps \ne 0 $.  

To motivate the approach that is to be followed, return to (\ref{eq:diracde}).
Except on the jump set $ {\mathcal J} := \left\{ t_1, t_2, \cdots, t_n \right\} $ , (\ref{eq:diracde}) would evolve smoothly according to the standard
ordinary differential equation (\ref{eq:unpde}).    In the ``impulsive differential equations'' viewpoint, a jump will occur at each
value $ t_i $, and this is usually {\em specified}  \cite{bainovsimeonov,barreiravalls,barreira,pan,hanlu,linyoungshear,linyoungkick,fennersiegmund,liuballinger}.  Thus, the system to
be examined would be (\ref{eq:unpde}) plus the jump maps specified at the times $ t_i $.  The relationship of each jump map
to the function $ g_i $ would be hidden in this approach.   Here, the intention is to reveal this connection (as done in other
studies \cite{catlla,griffithswalborn,belleyvirgilio_duffing,belleyvirgilio_lienard,dubeaukarrakchou}), while explicitly seeking
stable/unstable manifolds.
To retain the effect of the $ g_i $s and still make sense of equations such as (\ref{eq:diracde}), Catlla et~al \cite{catlla} suggest
the ``$\delta $-sequence'' approach which they apply to a first-order linear equation \cite{catlla}.  This idea, applied to the
present context, would necessitate the identification of a  ``$\delta $-family'' of functions $ \delta_\ell(t) $ which in the `limit'
$ \ell \downarrow 0 $ approach the Dirac delta $ \delta(t) $.  One way to specify this is to define this family as
piecewise continuous functions $ \delta_\ell: \R \rightarrow \R $ for $ \ell > 0 $, which have the property that for continuous functions $ h:
\R \rightarrow \R^n $, 
\begin{equation}
\fl \lim_{\ell \downarrow 0} \int_{-\infty}^{\infty} h(\tau) \, \delta_\ell(\tau-t) \, \d \tau = h(t) \, .
\label{eq:deltafamily}
\end{equation}
Then, a natural interpretation of (\ref{eq:diracde}) would be to look for solutions $ x_\ell(t) $ which satisfy
\begin{equation}
\dot{x}_\ell(t) = f \left( x_\ell(t) \right) + \eps \sum_{i=1}^n g_i \left( x_\ell(t), t \right) \delta_\ell(t-t_i)
\label{eq:xell}
\end{equation}
and subsequently take the limit $ \ell \downarrow 0 $ (if it exists).   The time values at which this limit becomes
difficult are the $ t_i $s at which the impulses occur.  While a pleasing implicit expression for the jumps in the solutions
occurring at these values for {\em any}
choice of the $ \delta $-family is possible when $ x $ is one-dimensional (see Proposition~5.1 in \cite{catlla}),  this
separations-of-variables approach cannot be used in this two-dimensional situation.

\begin{hypothesis}[Properties of impulsive perturbation]
\label{hyp:gdirac}
The perturbation in (\ref{eq:xell}) is associated with the following properties:
\begin{itemize}
\item[(a)] Define the {\em jump set} $ {\mathcal J} := \left\{ t_1, t_2, \cdots, t_n \right\} $, where the $ t_i $ are an increasing
set of values in $ \R $; 
\item[(b)] For each $ t \in \R $, $ g_i  \left( \centerdot, t \right) \in {\mathrm{C}}^1 \left( \Omega \right) $, with both $ g_i $ and $ D g_i $ bounded on $ \Omega $;
\item[(c)] For each $ x \in \Omega $, and each $ i \in \left\{ 1, 2, \cdots , n \right\} $, $ g_i (x, \centerdot) \in {\mathrm{C}}^1 \left(
\R \right) $.
\end{itemize}
\end{hypothesis}

\begin{lemma}[Existence, uniqueness, smoothness and invertibility of jump map]
\label{lemma:diracunique}
Choose the $ \delta $-family
\begin{equation}
\delta_\ell(t) = \frac{\alpha}{\ell} \I_{[- \ell,0 )}(t) + \frac{1- \alpha}{\ell} \I_{[0,\ell ]}(t) \quad , \quad \alpha \in [0,1] \, , 
\label{eq:deltaalpha}
\end{equation}
where $ \I $ is the indicator function,
and let $ G_i: \Omega \rightarrow \Omega $ be the ``jump map'' which takes the point $ x(t_i^-) $ to $ x(t_i^+) $ in (\ref{eq:diracde}).
Then, for $ \left| \eps \right| $ sufficiently small, $ G_i $ exists as a unique diffeomorphism on $ \Omega $ for each $ i $.
\end{lemma}

\proof
Fix an $ i \in \left\{ 1, 2, \cdots, n \right\} $ and choose an interval $ T_i $ containing $ t_i $ but none of the other points
from $ {\mathcal J} $.  Choose $ \ell > 0 $ small enough such that $ \left[t_i-\ell,t_i+\ell
\right] \subset T_i $.  Then, the dynamics in $ T_i $ are, from (\ref{eq:xell}),
\[
\fl \dot{x}_\ell(t) = f \left( x_\ell(t) \right) + \eps g_i \left( x_\ell(t), t \right) \delta_\ell(t-t_i) \quad , \quad t \in T_i \, .
\]
Integrating this from $ t_i - \ell $ to $ t_i + \ell $ yields
\begin{eqnarray*}
& & \fl x_\ell(t_i+\ell) - x_\ell(t_i-\ell) - \int_{t_i-\ell}^{t_i+\ell} f \left( x_\ell(t) \right) \, \d t \\
& = & \eps \left[ \frac{\alpha}{\ell} \int_{t_i-\ell}^{t_i} 
g_i \left( x_\ell(t),t \right) \, \d t + \frac{1-\alpha}{\ell} \int_{t_i}^{t_i+\ell} g_i \left( x_\ell(t),t \right) \, \d t \right] \\
& = & \eps  \left[ \alpha \! \! \int_{-1}^{0} 
g_i \left( x_\ell(t_i \! + \! \ell \tau),t_i \! + \! \ell \tau \right) \d \tau + (1-\alpha) \! \! \int_{0}^{1} g_i 
\left( x_\ell(t_i \! + \! \ell \tau),t_i \! + \! \ell \tau \right) \d \tau \right] \, .
\end{eqnarray*}
Letting $ x(t) = \lim_{\ell \downarrow 0} x_\ell(t) $ and taking the limit $ \ell \downarrow 0 $ above gives
\begin{equation}
x(t_i^+) - x(t_i^-) = \eps \left[ \alpha g_i \left( x(t_i^-),t_i \right) + \left( 1 - \alpha \right) g_i \left( x(t_i^+),t_i \right) \right] \, , 
\label{eq:jumpasymmetric}
\end{equation}
by taking into account the smoothness of the functions $ f $ and $ g_i $.  The expression (\ref{eq:jumpasymmetric}) is akin to
the idea of ``matched asymptotics'' \cite{catlla} which specifies a condition for the jump.  Thus, the mapping from $ x(t_i^-) $ to
$ x(t_i^+) $, if expressed as $ G_i $, is  defined implicitly on $ \Omega $ by
\[
\fl G_i(x) - x - \eps \left[ \alpha g_i \left( x ,t_i \right) + (1-\alpha) g_i \left( G_i (x), t_i\right) \right] = 0 \, ,
\]
the requirement is to find $ y = G_i(x) $ satisfying
\[
\fl y - x - \eps \left[ \alpha g_i(x,t_i) + (1-\alpha) g_i(y,t_i) \right] = 0 \, .
\]
Note that when $ \eps = 0 $, a unique solution for $ y(x,\eps) $ is $ y = x $.  Now, given Hypothesis~\ref{hyp:gdirac},
the $ y $-derivative of the left-hand side above differs from the identity by terms of size $ \O{\eps} $.
Thus for small enough $ \left| \eps \right| $ the determinant of this derivative matrix will be bounded away from zero, and the implicit
function theorem establishes that  for any $ x_0 \in \Omega $,
there exists an open neighbourhood $ B(x_0) $, and also a small interval containing $ 0 $
(say, $ E $), such that for $ (x,\eps) \in B(x_0) \times E $, $ y $ can be solved uniquely
as a function of $ (x,\eps) $.  This moreover establishes
that $ y $ is as smooth in $ x $ as is $ g_i $.  Since this works for any $ x_0 \in \Omega $,
a global smooth solution $ y(x,\eps) $ exists on $ \Omega \times E $. The argument in backwards time is similar,
establishing the existence, uniqueness and smoothness of $ G_i^{-1} $, which thereby proves the invertibility of each
$ G_i $.
\proofend

The choice of $ \delta $-family given in (\ref{eq:deltaalpha}) 
 incorporates the most common regularisation of the Dirac delta when 
$ \alpha = 1/2 $.   This is of a symmetric rectangular pulse \cite[e.g.]{griffithswalborn}.  The case $ \alpha = 1 $ 
is a simple one, in which case the jump map in the forward time direction automatically exists \cite[e.g.]{dubeaukarrakchou,catlla},
as can be seen from (\ref{eq:jumpasymmetric}).  (If $ \alpha = 0 $, this is true in the backward time direction.)
If using a different $ \delta $-family (such as tent functions or Gaussians), the proof of existence of
the jump map become more tricky.  Indeed, a (one-dimensional) example by Catlla et~al \cite{catlla} (their equation~
(5.14)) indicates that the jump map $ G_i $ may not exist for a general choice of the function $ g_i $.  There are however
different conditions from those given in Lemma~\ref{lemma:diracunique} under which existence in certain classes of
state-dependent impulsive
systems can be established \cite{belleyvirgilio_duffing,belleyvirgilio_lienard}.

Whenever $ t $ is well-removed from the jump set $ {\mathcal J} $, (\ref{eq:xell}) indicates that $ x(t) $ would simply
evolve according to $ \dot{x} = f(x) $.  As $ t $ crosses values in the jump set, a jump as given by (\ref{eq:jumpasymmetric}) needs to
be applied.  These factors can be combined in representing the $ \ell \downarrow 0 $ limit of (\ref{eq:xell}) in terms of
an integral equation.  In stating this,  it is possible to dispense with the 
explicit $ t $-dependence by redefining each $ g_i $ by $
g_i(x,t_i) \rightarrow g_i(x) $,
an abuse of notation which has shall be followed henceforth.  This leads to the {\em integral equation}
\begin{equation}
\fl x(t) = x(\beta) + \! \int_\beta^t  \! \! \! f \left( x(\xi) \right) \d \xi
+ \eps \sum_{i=1}^n
u \left( \beta, t_i, t \right)\,
\left[ \alpha g_i \! \left( x(t_i^-) \right) + (1- \alpha) g_i \! \left( x(t_i^+) \right) \right] \, ,
\label{eq:dirac}
\end{equation}
where 
\begin{equation}
\fl u \left(  \beta, t_i, t \right) := \left\{ \begin{array}{cl} 1  & \quad {\mathrm{if~}} \beta < t_i < t  \\
-1 & \quad {\mathrm{if~}}  \beta > t_i > t \\
0 & \quad {\mathrm{if~else}} \end{array} \right. \, .
\label{eq:timiddle}
\end{equation}
The remainder of this article focusses on (\ref{eq:dirac}), which is one particular rationalisation of the conceptual
form (\ref{eq:diracde}).  Within this approach, it will be possible to establish expressions for the impulsive analogues
of stable and unstable manifolds.

\begin{center}
\begin{figure}[t]
\includegraphics[width=0.4 \textwidth]{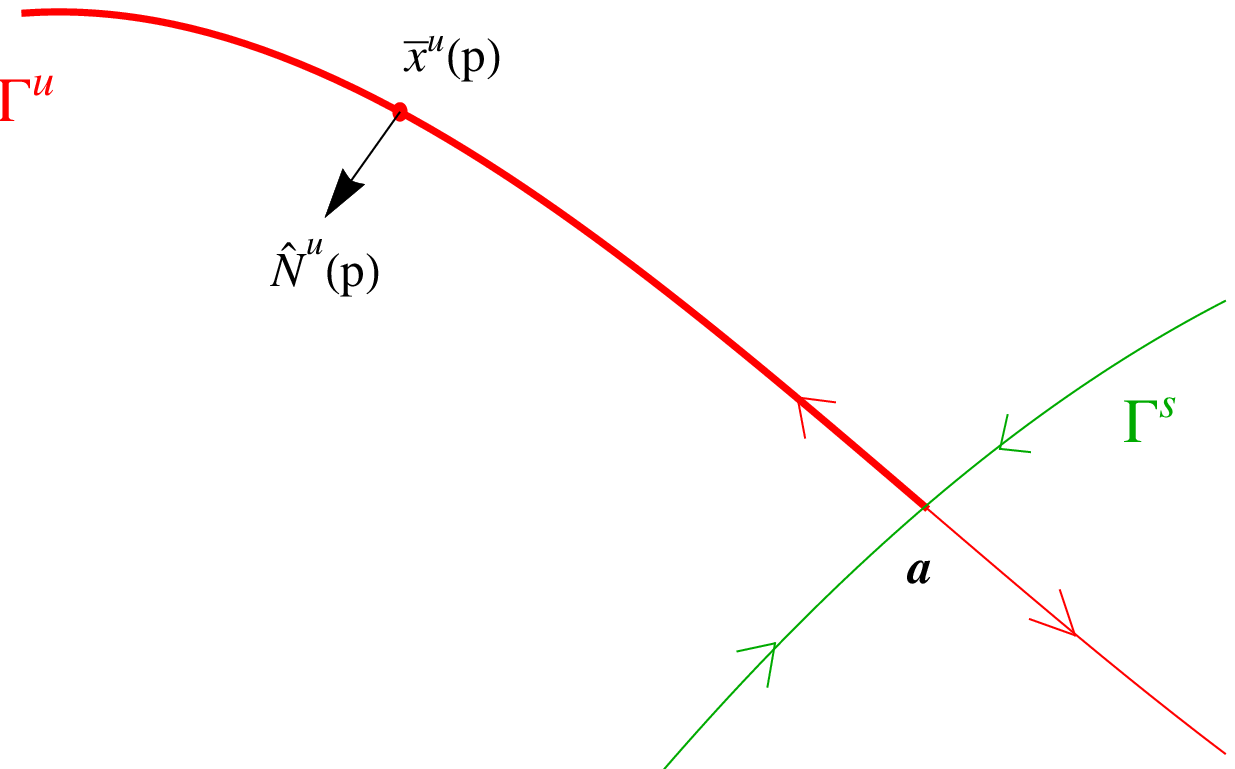}
\includegraphics[width=0.6 \textwidth]{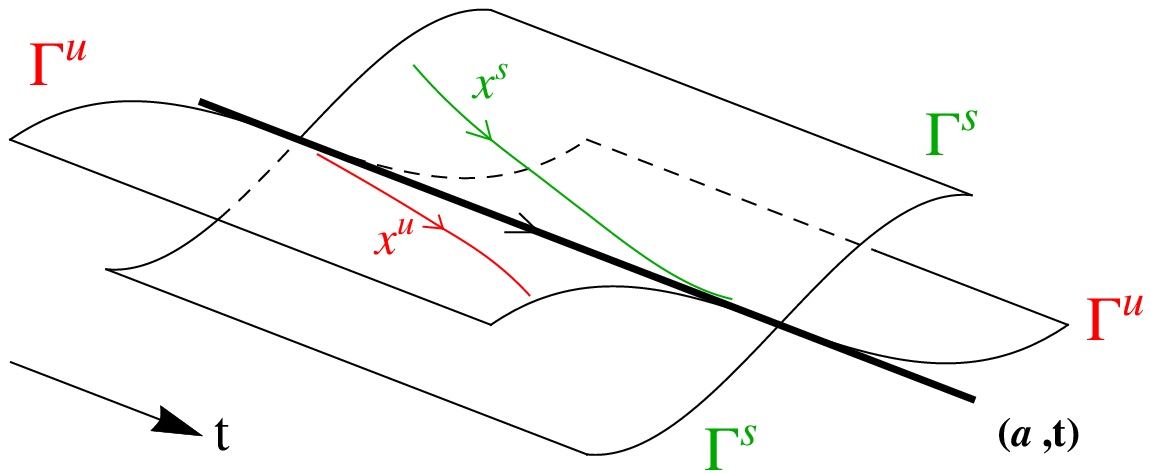}
\hspace*{0.2 \textwidth} (a) \hspace*{0.4 \textwidth} (b)
\caption{The $ \eps = 0 $ phase spaces for (\ref{eq:dirac}): (a) $ \Omega $, and (b) $ \Omega \times \R $, 
displaying hyperbolic trajectory [bold], and the two branches of each of the stable, $ \Gamma^s $,  and unstable, $ \Gamma^u $, manifolds.}
\label{fig:unp}
\end{figure}
\end{center}

\section{Pseudo-manifolds}
\label{sec:pseudomanifolds}

When $ \eps = 0 $, $ a $ was a saddle fixed point, with stable ($ \Gamma^s $) and unstable 
($ \Gamma^u $) manifolds existing as curves in $ \Omega $,
as shown in Figure~\ref{fig:unp}(a).   In the augmented $ (x,t) \in \Omega \times \R $ phase space this is representible as
a {\em hyperbolic trajectory} \cite{unsteady} $ (a,t) $ which possesses two-dimensional stable and unstable 
manifolds (also denoted by $ \Gamma^{s,u} $ with an abuse of notation), as shown in Figure~\ref{fig:unp}(b).  An important observation---to be useful later---is that the only points in
$ \Omega $ near $ a $ which in backwards time approach $ a $ are those lying on $ \Gamma^u $.
It is well-known that under smooth and bounded nonautonomous perturbations, $ (a,t) $ perturbs to 
$ \left( a_\eps(t), t \right) $, itself hyperbolic \cite{coppel,yi,tangential}.  In particular, this trajectory will retain its stable 
and unstable manifolds, for which it is possible to derive parametric expressions \cite{tangential}.

Under the nonsmooth integral equation evolution (\ref{eq:dirac}), however, the situation is different.  Consider choosing $ \beta < t_1 $ in (\ref{eq:dirac}), with $ x(\beta) = a $.  Since there is no perturbation to the steady flow until time $ t_1 $, 
 $ x(t) = a $ for $ t < t_1 $.  As $ t_1 $ is crossed, a jump to
$ x(t_1^+) = G_1 \left( x(t_1^-) \right) = G_1(a) $ will occur such that $ x(t_1^+) $ is $ \O{\eps}$-close to $ a $ but is generically
not the fixed point $ a $ of $ f $.  Thus, typically, the subsequent evolution of this trajectory will not be stationary. 
 At $ t_2 $, $ x(t) $ will once again jump, and so on, until passing the final jump time $ t_n $.  Since $ x(t_n^+) $ will
also not be a fixed point, the subsequent evolution will be governed by $ x(t) = x(t_n^+) + \int_{t_n}^t f \left( x(\tau) \right) \, \d
\tau $, and will generically experience exponential separation from $ a $ since $ a $ is unstable. 
This trajectory, labelled $ a_+(t) $, will be defined for $ t \in (-\infty, T_u] \setminus {\mathcal J} $ for any 
finite $ T_u $ as long as the trajectory remains within $ \Omega $, and be $ \O{\eps} $-close to $ a $ in this domain of validity.  In a 
similar vein,  $ a_-(t) $ will be the trajectory obtained
by taking $ x(\beta) = a $ for $ \beta > t_n $, and evolving (\ref{eq:dirac}) {\em backwards} in time; this will be defined for
$ t \in [T_s, \infty) \setminus {\mathcal J} $ for $ -T_s $ arbitrarily large but finite.  The two trajectories $ a_+(t) $ and $ a_-(t) $
are respectively $ a $'s forwards and backwards iterates under the perturbed flow, and will not coincide in general; the unique
hyperbolic trajectory $ a_\eps(t) $ present in the smooth situation does not occur.  Therefore, stable and unstable manifolds
attached to $ \left( a_\eps(t), t \right) $ in the standard nonautonomous sense cannot be defined, and indeed the lack of continuity of the $ a_\pm(t) $ trajectories questions the very usage of the term `manifolds.'

\begin{definition}[Unstable pseudo-manifold]
\label{definition:pseudomanifoldunstable}
The unstable pseudo-manifold of $ a $ 
in the augmented phase space $ \Omega \times (-\infty, T_u] \setminus {\mathcal J} $ for any finite $ T_u $ is defined by 
\begin{equation}
\fl \Gamma_\eps^u  := \bigcup_{\beta \in (-\infty, T_u] \setminus {\mathcal J}} \left\{ \left(x(\beta), \beta \right) ~:~
{\mathrm{all}}~x(\beta) \in \Omega~{\mathrm{for~which}}~x(t) \rightarrow a~{\mathrm{as}}~t \rightarrow - \infty  \right\} \, , 
\label{eq:pseudomanifoldunstable}
\end{equation}
where $ x(t) $ is the evolution defined in (\ref{eq:dirac}).
\end{definition}

The explanation for Definition~\ref{definition:pseudomanifoldunstable} appears in Figure~\ref{fig:unstable}, where in the diagram,
only one of the two branches (that corresponding to the upper left surface in Figure~\ref{fig:unp}) is shown.  Only the 
first two jump values, at $ t = t_1 $ and $ t_2 $, are displayed.  All points on $ \Gamma_\eps^u $ for $ t < t_1 $ will decay in
backwards time to $ a $, since in this region the situation is exactly as in Figure~\ref{fig:unp}, with $ \Gamma_\eps^u $ 
coinciding with $ \Gamma^u $.  At $ t = t_1^- $, $ \Gamma_\eps^u $ forms a curve in the time-slice $ t = t_1 $.  However,
all points on this curve jump according to the map $ G_1 $ because of the impulse, thereby forming a new curve (that
corresponding to $ t = t_1^+ $) in the time-slice $ t_1 $.  These points then evolve continuously according to the vector field
$ f $ until $ t_2 $, whereupon $ G_2 $ applies, to create another curve.  Since it is only points on the collection of
surfaces $ \Gamma_\eps^u $ which get mapped back to $ \Gamma_u $ for $ t < t_1 $, it is exactly points on these surfaces
which attracted towards $ a $ in backwards time.  It should be noted that the `special' trajectory $ \left( a_+(t), t \right) $
is a boundary of $ \Gamma_\eps^u $.  It is in fact a `hyperbolic-like trajectory' in backwards time {\em only}, in the sense
that points on the attached surface get attracted towards it at an exponential rate in backwards time.
Given the discontinuities, $ \Gamma_\eps^u $ fails to be well-defined on
the time-slices $ t = t_i $, $ t_i \in {\mathcal J} $.  This lack of smoothness of $ \Gamma_\eps^u $ is what prompts the
term {\em pseudo}-manifold in Definition~\ref{definition:pseudomanifoldunstable}.   A similar definition is therefore possible
for the stable pseudo-manifold:

\begin{figure}
\begin{center}
\includegraphics[width=\linewidth]{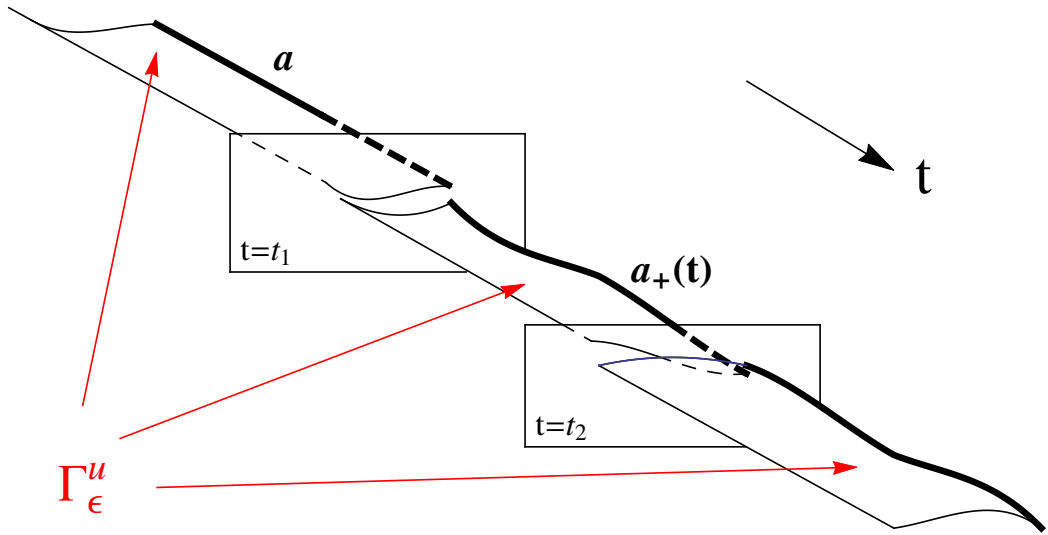} 
\caption{The unstable pseudo-manifold $ \Gamma_\eps^u $ of $ a $ associated with (\ref{eq:dirac}), which comprises
segments of smooth surfaces with jump discontinuities at $ t_i $, $ i = 1, 2, \cdots, n $.  The thick curve is the `hyperbolic-like'
trajectory $ a_+(t) $, to which trajectories on $ \Gamma_\eps^u $ are attracted in backwards time.}
\label{fig:unstable}
\end{center}
\end{figure}

\begin{definition}[Stable pseudo-manifold]
\label{definition:pseudomanifoldstable}
The stable pseudo-manifold of $ a $ 
in the augmented phase space $ \Omega \times [T_s,\infty) \setminus {\mathcal J} $ is defined by 
\begin{equation}
\fl \Gamma_\eps^s  := \bigcup_{\beta \in [T_s,\infty) \setminus {\mathcal J}} \left\{ \left(x(\beta), \beta \right) ~:~
{\mathrm{all}}~x(\beta) \in \Omega~{\mathrm{for~which}}~x(t) \rightarrow a~{\mathrm{as}}~t \rightarrow \infty  \right\} \, , 
\label{eq:pseudomanifoldstable}
\end{equation}
where $ x(t) $ is the evolution defined in (\ref{eq:dirac}).
\end{definition}

Now, an expression for the pseudo-manifolds is sought.  This shall be expressed in a parametric way, and the unstable
pseudo-manifold shall be the initial focus.  When $ \eps = 0 $, the unstable manifold might be thought of in terms of a solution
 $ \bar{x}^u(t) $ to (\ref{eq:unpde}) which satisfies $ \bar{x}^u(t) \rightarrow a $ as $ t \rightarrow - \infty $.  Thus, $ p \in
(-\infty, P] $, for $ P $ as large as desired but finite, can be used to parametrise a segment of the unstable manifold in $ \Omega $
in the form $ \bar{x}^u(p) $, as indicated in Figure~\ref{fig:unp}(a).  Since a finiteness assumption on $ P $ is imposed, this means that this situation captures varied
possibilities for the `other end of the manifold,' which might attach to another fixed point, escape to infinity, spiral in towards a 
limit cycle, etc.  By having $ P $ finite, the unstable manifold is clipped at some point; the curve of interest, in $ \Omega $, is
of finite length.  When considering this in the augmented $ \Omega \times \R $ phase space, one
would have trajectories $ \left( \bar{x}^u(t), t \right) $ lying on the two-dimensional unstable manifold, of which one is shown
in Figure~\ref{fig:unp}(b).  Indeed, {\em all} the trajectories on this manifold can be obtained by simply shifting this one 
trajectory, since the system is autonomous.  Another way to think of this is that for each initial condition chosen on $ \Gamma^u $
in the $ \Omega $ phase-space of Figure~\ref{fig:unp}(a) will generate a trajectory on the unstable manifold in Figure~\ref{fig:unp}(b).

Now suppose $ \eps \ne 0 $.  Consider a fixed time slice $ t 
\in (-\infty, T_u] $
on the augmented $ \Omega \times \R \setminus {\mathcal J} $ phase-space.  Within this time-slice, if $ \eps = 0 $, the
picture of the unstable manifold would be as shown in Figure~\ref{fig:unp}(a).  Thus, $ p \in (-\infty,P] $ will characterise a location
on the unperturbed unstable manifold.  Of course, after perturbation, the unstable pseudo-manifold will not lie exactly on
$ \Gamma^u $.  At the point $ \bar{x}^u(p) $ (i.e., the $ p $-parametrisation point), consider drawing a normal to the unperturbed
unstable manifold in the direction given by $ f^\perp \left( \bar{x}^u(p) \right) $.  Here, the perpendicular notation indicates a
rotation of a two-dimensional vector by $ \pi/2 $ in the anticlockwise direction, and since $ f $ is parallel to the original manfold,
$ f^\perp $ shall be normal to it.  More specifically, referring to Figure~\ref{fig:unp}(a), define
\begin{equation}
\fl \hat{N}^u(p) := \frac{ f^\perp \left( \bar{x}^u(p) \right) }{ \left| f \left( \bar{x}^u(p) \right) \right|} \quad , \quad f^\perp :=
\left( \begin{array}{cc} 0 & -1 \\ 1 & 0 \end{array} \right) f \, . 
\label{eq:normalu}
\end{equation}

\begin{theorem}[Unstable pseudo-manifold]
\label{theorem:unstabledirac}
Consider (\ref{eq:dirac}) under Hypotheses~\ref{hyp:fdirac} and \ref{hyp:gdirac}.  The unstable pseudo-manifold
of $ a $ has a parametric representation $ \left( x_\eps^u(p,t), t \right) $ 
with parameters $ (p,t) \in (-\infty, P] \times (-\infty,T_u] \setminus {\mathcal J} $ for
arbitrarily large but fixed $ P $ and $ T_u $, such that 
\begin{equation}
\fl \left[ x_\eps^u(p,t) - \bar{x}^u(p) \right] \cdot \hat{N}^u(p)  = \eps \frac{M^u(p,t)}{\left| 
f \left( \bar{x}^u(p) \right) \right|} + \O{\eps^2}\, , 
\label{eq:unstabledirac}
\end{equation}
where the associated unstable Melnikov function is given by
\begin{equation}
\fl M^u(p,t) = \sum_{i=1}^n \I_{(t_i,\infty)}(t) j_i^u(p,t) + \sum_{i=1}^{{\mathrm{max}} \left\{ j: t_j < t \right\}} \int_{t_i}^t R_p^u(t-\xi) j_i^u(p,\xi) \, \d \xi \, , 
\label{eq:mudirac}
\end{equation}
in which
\begin{equation}
\fl j_i^u(p,t) :=  f^\perp \left( \bar{x}^u(t_i-t+p) \right) \cdot g_i \left(  \bar{x}^u(t_i-t+p)  \right) 
\label{eq:jumpu}
\end{equation}
and the resolvent $ R_p^u $ is defined in terms of Laplace transforms with respect to $ t $ by
\begin{equation}
\fl R_p^u(t) := {\mathcal L}^{-1} \left\{ \frac{\hat{F}_p^u(s)}{1 - \hat{F}_p^u(s)} \right\}(t) \quad , \quad
\hat{F}_p^u(s) := {\mathcal L} \left\{ {\mathrm{Tr}} \, D f \left( \bar{x}^u(p-t) \right) \right\}(s) \, .
\label{eq:resolventu}
\end{equation} 
\end{theorem}

\proof
See \ref{sec:melnikov_impulsive_unstable}.
\proofend

\begin{remark}[Independence on asymmetry of Dirac impulse formulation] \normalfont
An interesting feature of the leading-order normal displacement of the unstable manifold, as given in Theorem~\ref{theorem:unstabledirac}, is that it is independent of $ \alpha $.  Thus asymmetric interpretations of a Dirac
impulse (in the form of (\ref{eq:dirac})) do not affect this quantity.  It is likely that the higher-order terms 
in the displacement are, however, dependent on $ \alpha $.
\end{remark}

\begin{corollary}[Unstable pseudo-manifold under area-preservation]
\label{corollary:unstabledirac}
Under the conditions of Theorem~\ref{theorem:unstabledirac}, consider the additional assumption that $ f $ is 
area-preserving.  Then, (\ref{eq:mudirac}) simplifies to
\begin{equation}
\fl M^u(p,t) = \sum_{i=1}^n \I_{(t_i,\infty)}(t) f^\perp \left( \bar{x}^u(t_i-t+p) \right) \cdot g_i \left(  \bar{x}^u(t_i-t+p) \right) \, .
\label{eq:mudiracareapreserving}
\end{equation}
\end{corollary}

\proof
Since in this case $ {\mathrm{Tr}} \, D f = 0 $, $ R_p^u = 0 $ from (\ref{eq:resolventu}).  Thus, from (\ref{eq:mudirac}), 
$ M^u(p,t) = j_p^u(t) $ directly.
\proofend

\begin{remark}[Formal Melnikov computation under impulses] \normalfont
\label{remark:impulseformal}
If the situation being considered is $ \dot{x} = f(x) + \eps g(x,t) $ where $ g(x,t) $ is smooth, then the distance 
expression (\ref{eq:unstabledirac}) for the normal displacement of the unstable manifold continues to hold, but now with
\begin{equation}
\fl M^u(p,t) =  \int_{-\infty}^t \! \! \exp \left[ \int_{\xi-t+p}^p {\mathrm{Tr}} \left[  D f \left( \bar{x}^u(\tau) \right) \right] \d \tau \right] f^\perp \left( 
\bar{x}^u(\xi \! - \! t \! + \! p) \right) \cdot g \left( \bar{x}^u(\xi\! -\! t \! + \!p), \xi \right) \, \d \xi
\label{eq:musmooth}
\end{equation}
as shown in \cite{tangential}.  If $ {\mathrm{Tr}} \, D f = 0 $, then a purely formal replacement of $ g(x,t) $ above with
$ \sum_{i=1}^n \delta(t-t_i) g_i(x) $ directly gives the formula (\ref{eq:mudiracareapreserving}).  It is however instructive that
the same formal approach gives the {\em wrong} result (i.e., not (\ref{eq:mudirac})) if $ {\mathrm{Tr}} \, D f \ne 0 $, thereby
highlighting the necessity of following the integral equation approach.
\end{remark}

The modifications for the {\em stable} pseudo-manifold
are analogous.  When $ \eps = 0 $, $ \bar{x}^s(t) $ is assumed to be a trajectory on a branch of the stable manifold, such that
$ \bar{x}^s(t) \rightarrow a $ as $ t \rightarrow \infty $.  The normal vector shall be defined by $ \hat{N}^s(p) := f^\perp \left(
\bar{x}^s(p) \right) / \left| f \left( \bar{x}^s(p) \right) \right| $.
When $ \eps \ne 0 $, $ a $ remains a fixed point for $ t > t_n $, with its stable manifold well-defined; this is simply
taken in backwards time across the time-discontinuities to generate the stable pseudo-manifold.  The leading-order
representation of its normal displacement, just as for the unstable pseudo-manifold, is independent of $ \alpha $:

\begin{theorem}[Stable pseudo-manifold]
\label{theorem:stabledirac}
Consider (\ref{eq:dirac}) under Hypotheses~\ref{hyp:fdirac} and \ref{hyp:gdirac}.  The stable pseudo-manifold
of $ a $ has a parametric representation $ \left( x_\eps^s(p,t), t \right) $ 
with parameters $ (p,t) \in [-P,\infty) \times [T_s,\infty) \setminus {\mathcal J} $ for
arbitrarily large but fixed $ P $ and $ -T_s $, such that 
\begin{equation}
\fl \left[ x_\eps^s(p,t) - \bar{x}^s(p) \right] \cdot \hat{N}^s(p)  = \eps \frac{M^s(p,t)}{\left| 
f \left( \bar{x}^s(p) \right) \right|} + \O{\eps^2}\, , 
\label{eq:stabledirac}
\end{equation}
where the associated stable Melnikov function is given by
\begin{equation}
\fl M^s(p,t) = - \sum_i^n \I_{(-\infty,t_i)}(t) j_i^s(p,t) + \sum_{i={\mathrm{min}} \left\{ j: t_j > t \right\}}^n 
\int_t^{t_i} R_p^s(t-\xi) j_i^s(p,\xi) \, \d \xi \, , 
\label{eq:msdirac}
\end{equation}
in which 
\begin{equation}
\fl j_i^s(p,t) :=  f^\perp \left( \bar{x}^s(t_i-t+p) \right) \cdot g_i \left(  \bar{x}^s(t_i-t+p) \right) 
\label{eq:jumps}
\end{equation}
and the resolvent $ R_p^s $ is defined by
\begin{equation}
\fl R_p^s(t) := {\mathcal L}^{-1} \left\{ \frac{\hat{F}_p^s(s)}{1 + \hat{F}_p^s(s)} \right\}(-t) \quad , \quad
\hat{F}_p^s(s) := {\mathcal L} \left\{ {\mathrm{Tr}} \, D f \left( \bar{x}^s(p+t) \right) \right\}(s) \, .
\label{eq:resolvents}
\end{equation} 
\end{theorem}

\proof
While this is in principle similar to Theorem~\ref{theorem:unstabledirac}, the fact that the functions are defined on $ \R^- $ as opposed to $ \R^+ $ require subtle adjustments when using the Laplace transform;  details are outlined in\ref{sec:melnikov_impusive_stable}.
\proofend

\begin{corollary}[Stable pseudo-manifold under area-preservation]
\label{corollary:stabledirac}
Under the conditions of Theorem~\ref{theorem:stabledirac}, consider the additional assumption that $ f $ is 
area-preserving.  Then, (\ref{eq:msdirac}) simplifies to
\begin{equation}
\fl M^s(p,t) = - \sum_{i=1}^n \I_{(-\infty,t_i)}(t) f^\perp \left( \bar{x}^s(t_i-t+p) \right) \cdot g_i \left(  \bar{x}^s(t_i-t+p) \right) \, .
\label{eq:msdiracareapreserving}
\end{equation}
\end{corollary}

\proof
Simply set $ {\mathrm{Tr}} \, D f = 0 $, as in the proof of Corollary~\ref{corollary:unstabledirac}.
\proofend

\begin{remark}[Pseudo-manifolds and unsteady transport barriers] \normalfont
\label{remark:flowbarriers}
Stable and unstable manifolds in the unsteady infinite-time context form transport barriers 
in unsteady flows \cite{peacockhaller,siam_book,eigenvector,aperiodic,unsteady}; an interpretation of this will be provided in Section~\ref{sec:flux}.
There is considerable ongoing work in determining analogous entities in time-dependent flows which are known only
over a {\em finite-time}, in which finite-time
versions of properties associated with stable/unstable manifolds are used to determine these barriers.  For example,
the exponential attraction/repulsion property is captured in seeking finite-time Lyapunov exponents \cite{shadden}; curves/surfaces
of extremal attraction/repulsion in the definition of hyperbolic Lagrangian coherent structures \cite{hallerreview}; flow 
separating property in transfer operator methods \cite{froylandpadberg}; tangent vectors to manifolds associated with
Oseledets splitting \cite{eigenvector}; etc.  It is not clear whether these different diagnostic approaches for
determining flow barriers would be practicable in instances in which the system had impulses; however, the pseudo-manifold definitions
given here do indeed enjoy the same the transport barrier properties that are associated with standard stable/unstable manifolds.
\end{remark}

\section{Example: parabolic pseudo-manifolds}
\label{sec:parabolic}

\begin{figure}
\begin{center}
\includegraphics[width=0.47 \linewidth]{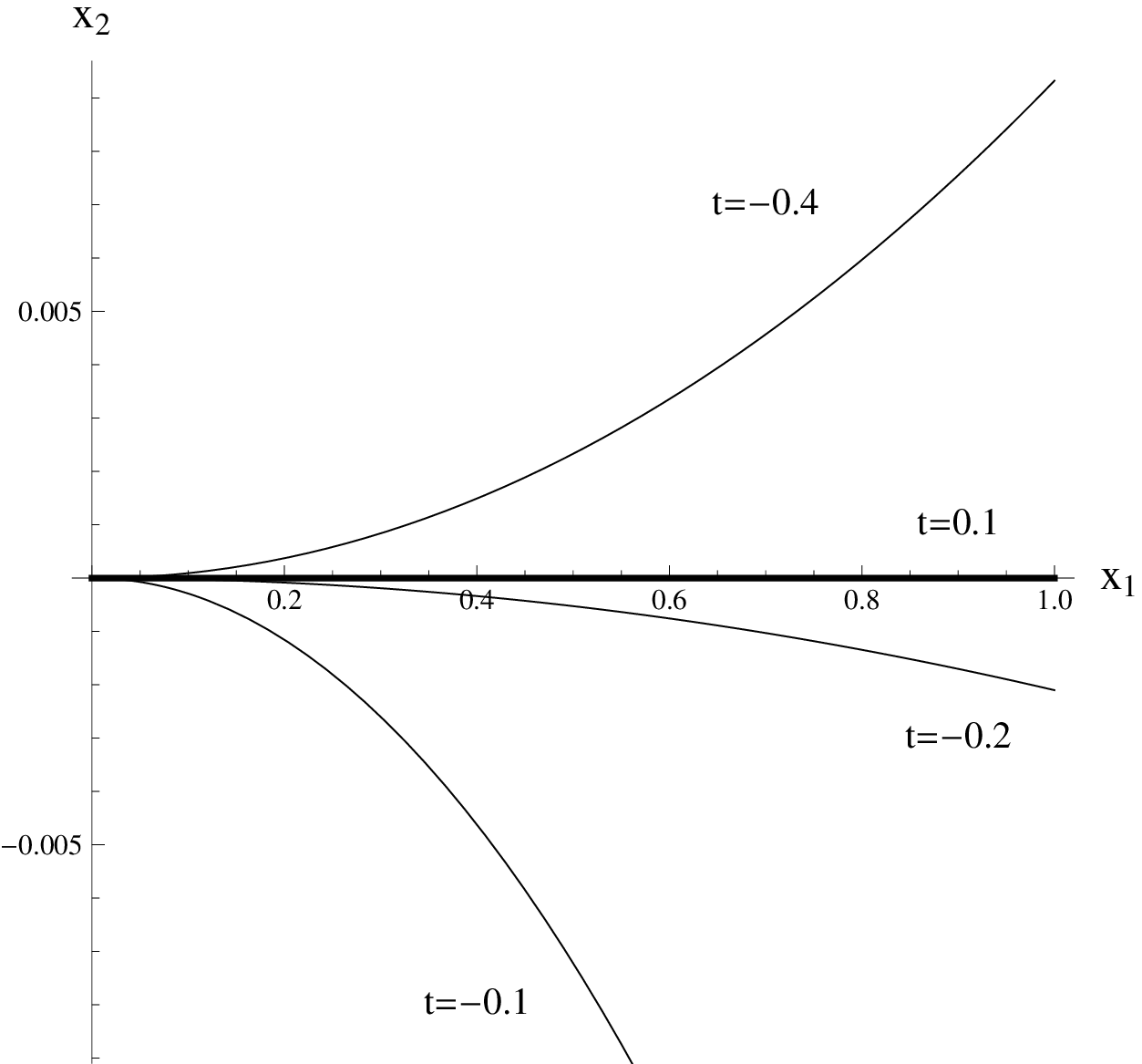} 
\includegraphics[width=0.47 \linewidth]{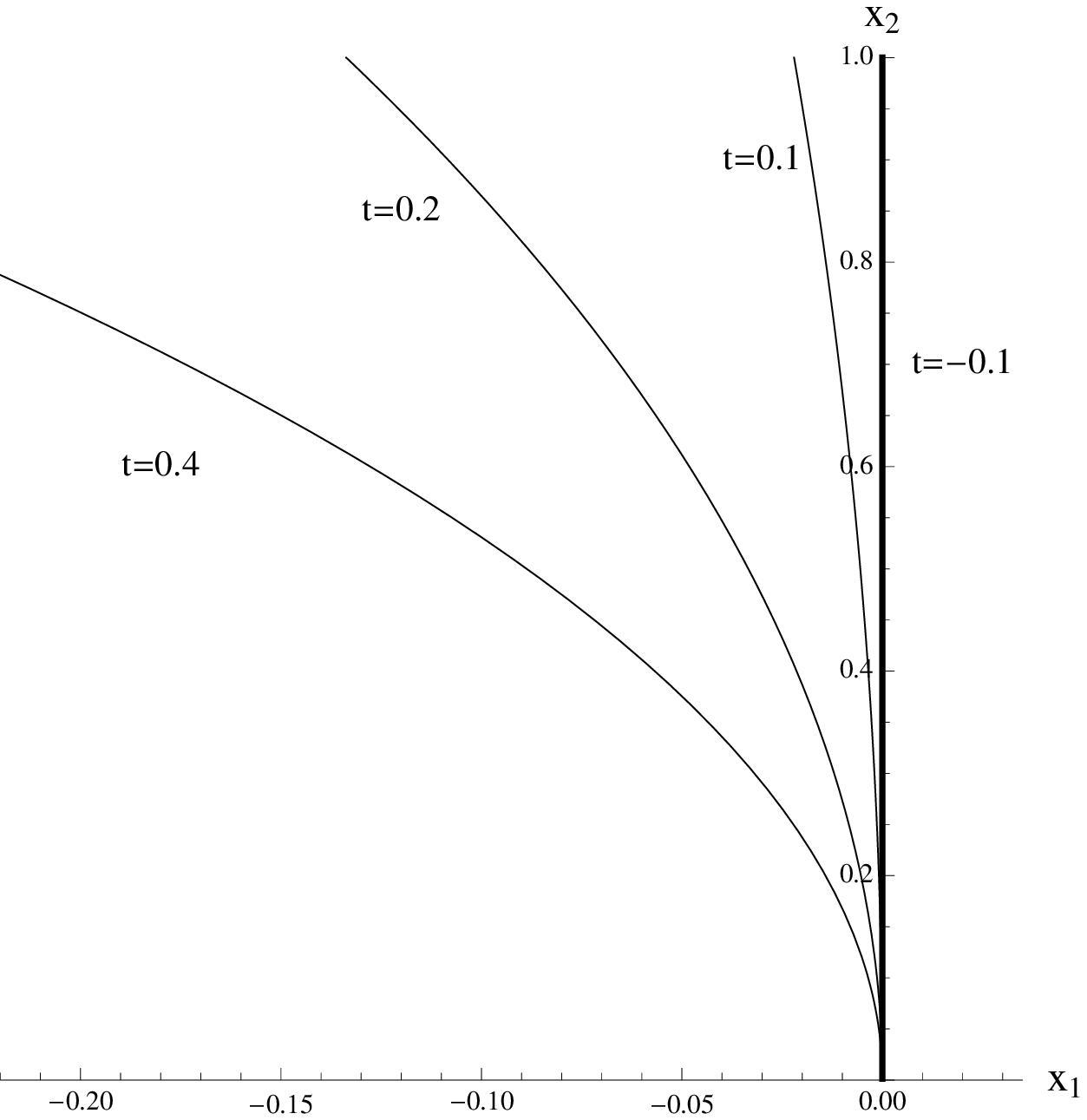} 
\caption{The stable pseudo-manifold (left) and unstable pseudo-manifold (right) 
for the example in Section~\ref{sec:parabolic} with $ \eps = 0.1 $,  at different $ t $-values.}
\label{fig:parabolicimpulse}
\end{center}
\end{figure}

Suppose $ f(x) = (-3 x_1, x_2) $, which corresponds to a saddle point at the origin with stable and unstable manifolds along
the $ x_1 $ and $ x_2 $ axes.   For the branch of the stable manifold lying along the positive $ x_1 $ axis,
\[
\fl \bar{x}^s(t) = \left( \begin{array}{c} e^{-3 t} \\ 0 \end{array} \right) \, , \, 
f \left( \bar{x}^s(t) \right) =  \left( \begin{array}{c} -3 e^{-3 t} \\ 0 \end{array} \right) \, , \, 
f^\perp \left( \bar{x}^s(t) \right) = \left( \begin{array}{c} 0 \\ -3 e^{-3 t} \end{array} \right) \, .
\]
For simplicity, suppose there is only one impulse occurring at $ t_1 = 0 $, with corresponding $ g_1(x,t)  = 
\left( x_1^2 + x_2^2, x_1^2 \cos t \right) $.  Then, 
\[
\fl j_1^s(p,t) =  f^\perp \cdot g_1 \left( \bar{x}^s(t_1 - t + p), t_1 \right) 
=   -3 e^{-3(0-t+p)} e^{-6(0-t+p)}  \cos 0 = -3  e^{-9 p}  e^{9 t} 
\]
For this $ f $,  $ {\mathrm{Tr}} \, D f = -3 +1 = -2 $.  Thus, $ \hat{F}_p^s(s) = -2/s $, and 
\[
\fl R_p^s(t) = {\mathcal L}^{-1} \left\{ \frac{-2/s}{1+(-2)/s} \right\}(t) = -2 e^{2 t} \, . 
\]
Using (\ref{eq:msdirac}),
\begin{eqnarray*}
\fl M^s(p,t) & = & - \I_{(-\infty,0)}(t) (-3 e^{-9 p} e^{9 t})+  \I_{(-\infty,0)}(t) \int_t^0 (-2 e^{2(t-\xi)}) (-3   e^{-9 p}  e^{9 \xi}) \, \d \xi \\
& = &  3 e^{-9 p} \I_{(-\infty,0)}(t) \left(  e^{9 t} + 2 e^{2 t} \frac{e^{7 \xi}}{7} \Big]_{\xi=0}^t \right) \\
& = & 3 e^{-9 p} \I_{(-\infty,0)}(t) \left( \frac{9}{7} e^{9 t} - \frac{2}{7} e^{2 t} \right) \, .
\end{eqnarray*}
The component of the stable manifold is, from (\ref{eq:stabledirac}), 
\begin{eqnarray*}
\fl \left[ x_\eps^s(p,t) - \bar{x}^s(p) \right] \cdot \left( -  \hat{x}_2 \right) & = & \eps \frac{ 3 e^{-9 p} \I_{(-\infty,0)}(t) \left\{ \frac{9}{7} e^{9 t} - \frac{2}{7} e^{2 t} \right\} }{ 3 e^{-3  p}} + \O{\eps^2} \\
& = & \eps \frac{e^{-6 p}}{7} \I_{(-\infty,0)}(t) \left[ 9 e^{9 t} -2  e^{2 t} \right]  + \O{\eps^2} \, .
\end{eqnarray*}
This means that the stable pseudo-manifold is
\[
\fl \tilde{\Gamma}_\eps^s = \left\{ \left( \left( \begin{array}{c} e^{-3 p} + \O{\eps} \\ \eps \frac{e^{-6 p}}{7} \I_{(-\infty,0)}(t) \left[ 2 e^{2 t} - 9 e^{9 t} \right]  + \O{\eps^2}
\end{array} \right), t \right) ~: ~t > T_s ~,~p > -P  \right\} \, ,
\]
where the $ \O{\eps} $ term in the $ x_1 $-component is since the theory only manages to capture the $ \O{\eps} $-normal
component of the manifold displacement; in general, there will also be a $ \O{\eps} $ modification in the tangential
direction (which has been quantified for smooth perturbations \cite{tangential}).
While the above is a $ (p,t) $ parametrisation for the stable pseudo-manifold (to leading-order), a formula for
the stable pseudo-manifold curves in each time-slice is easily obtained by eliminating $ p $ from the above, which gives
\[
\fl x_2 = \eps \frac{x_1^2}{7} \I_{(-\infty,0)}(t) \left[ 2 e^{2 t} - 9 e^{9 t} \right] + \O{\eps^2} \quad , \quad (x_1 > 0) \, .
\]
It is apparent that the tangential component becomes irrelevant to leading-order in this formulation.
Thus, while the stable pseudo-manifold is a straight line along the $ x_1 $ axis for $ t > 0 $, as $ t $ crosses $ 0 $ it abruptly
switches to approximately a parabolic curve initially given by $ x_2 = - \eps x_1^2 $ for $ x_1 > 0 $.  As time becomes additionally negative, 
the curvature of this parabolic curve evolves, as shown in the left panel of Figure~\ref{fig:parabolicimpulse}.  It is interesting
to note that the coefficient of the parabolic term changes sign when $ t = ( \ln 2/9 )/7 $ (approximately $ - 0.215 $), which means
that the parabola which opened `downwards' for negative $ t $-values near $ 0 $, opens `upwards' for more negative values.  Now what is important about these curves is if conditions were chosen on them at the labelled time, their trajectories
will eventually approach the origin as $ t \rightarrow \infty $.

Next, the unstable pseudo-manifold which perturbs from the unstable manifold branch lying along the $ + x_2 $ axis
is considered.  In this case, 
\[
\fl \bar{x}^u(t) = \left( \begin{array}{c} 0 \\ e^{t} \end{array} \right) \, , \, 
f \left( \bar{x}^u(t) \right) =  \left( \begin{array}{c} 0 \\ e^{t} \end{array} \right) \, , \, 
f^\perp \left( \bar{x}^u(t) \right) = \left( \begin{array}{c} - e^{t} \\ 0 \end{array} \right) \, , 
\]
and from (\ref{eq:jumpu}), 
\[
\fl j_1^u(p,t) =   f^\perp \cdot g_1 \left( \bar{x}^u(t_1 - t + p), t_1 \right) 
=  - e^{0-t+p} e^{2(0-t+p)} = - e^{3 p} e^{-3 t} \, .
\]
The relevant resolvent is, from (\ref{eq:resolventu}), 
\[
\fl R_p^u(t) = {\mathcal L}^{-1} \left\{ \frac{-2/s}{1-(-2)/s} \right\}(-t) = - 2 e^{-2 t} \, , 
\]
from which, using (\ref{eq:mudirac}), 
\begin{eqnarray*}
\fl M^u(p,t) & = &  \I_{(0,\infty)}(t) \left( - e^{3 p} e^{-3 t} \right) +  \I_{(0,\infty)}(t) \int_0^t (-2 e^{-2 (t-\xi)} ) (  - e^{3 p} e^{-3 \xi}  )  \, \d \xi \\
& = &  \I_{(0,\infty)}(t) e^{3p} \left[ 2 e^{2 t} - 3 e^{-3 t} \right]  \, . 
\end{eqnarray*}
The unstable pseudo-manifold expression (\ref{eq:unstabledirac}) therefore gives
\[
\fl \left[ x_\eps^u(p,t) - \bar{x}^u(p) \right] \cdot (-\hat{x}_1) =   \eps \frac{e^{3 p}}{e^p} \I_{(0,\infty)}(t) \left[ 2e^{2 t} - 3 e^{-3 t} \right]  + \O{\eps^2} \, .
\]
The $ \O{\eps} $ parametric approximation for the stable pseudo-manifold is therefore
\[
\fl \tilde{\Gamma}_\eps^u = \left\{ \left( \left( \begin{array}{c} \eps e^{2 p} \I_{(0,\infty)}(t) \left[ 3 e^{-3  t} - 2 e^{ 2 t} \right] + \O{\eps^2}
 \\ e^{p} + \O{\eps} \end{array} \right), t \right) ~: ~t < T_u ~,~p < P  \right\} \, ,
\]
and the nonparametric form is
\[
\fl x_1 = \eps x_2^2  \I_{(0,\infty)}(t) \left[ 3 e^{-3  t} - 2 e^{ 2 t} \right] + \O{\eps^2} \quad , \quad (x_2 > 0) \, , 
\]
which is also parabolic to leading-order, but now for $ t > 0 $. This is shown in the right panel of Figure~\ref{fig:parabolicimpulse}.

\section{Persistent heteroclinic trajectories}
\label{sec:heteroclinic}

Consider again (\ref{eq:dirac}) under Hypotheses~\ref{hyp:fdirac} and \ref{hyp:gdirac}, with the following additional hypothesis:

\begin{hypothesis}[Heteroclinic connection]
\label{hyp:heteroclinic}
The unperturbed ($\eps = 0 $) system (\ref{eq:dirac}) also satisfies 
\begin{itemize}
\item[(a)] There exists $ b \in \Omega $ (which might be the same point as $ a $) such that $ f(b) = 0 $ and $ D f(b) $ possesses
a positive and a negative eigenvalue;
\item[(b)] When considered in the $ \Omega $ phase-space, a branch of the unstable manifold of $ a $ coincides with a branch of the stable manifold of $ b $, forming a {\em heteroclinic manifold} $ \Gamma $ which can be parametrised by $ \bar{x}(p) $, $ p \in \R $ such that $ \bar{x}(p) \rightarrow a $ as $ p \rightarrow - \infty $
and $ \bar{x}(p) \rightarrow b $ as $ p \rightarrow \infty $, where $ \bar{x}(t) $ is a solution to (\ref{eq:dirac}) when $ \eps = 0 $.
\end{itemize}
\end{hypothesis}

\begin{figure}[t]
\centering
\includegraphics[width=0.8 \textwidth]{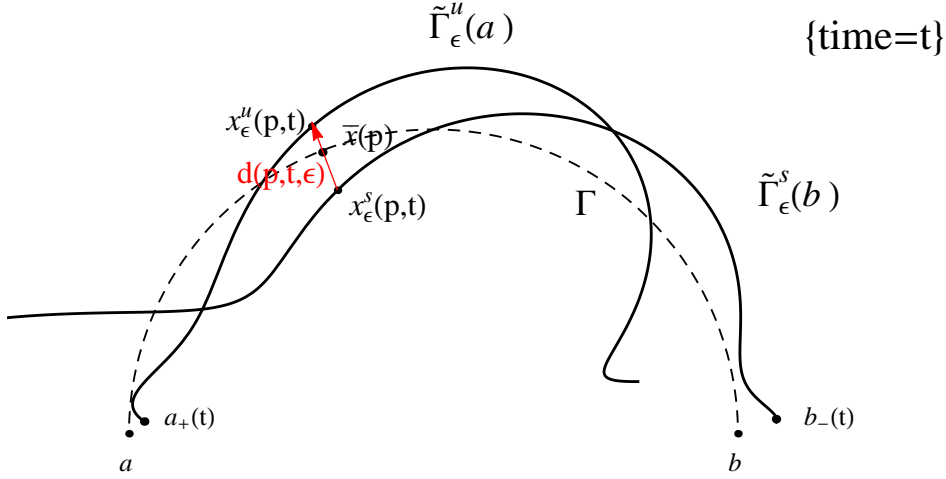}
\caption{Unperturbed heteroclinic manifold $ \Gamma $ [dashed] with the perturbed pseudo-manifolds in the time-slice $ t $;
the signed distance $ d(p,t,\eps) $ measured in the $ f^\perp $ normal direction to $ \Gamma $ at $ \bar{x}(p) $ is specified in 
Theorem~\ref{theorem:heteroclinic}.}
\label{fig:distance}
\end{figure}

If the point $ b $ is the same as $ a $, what is being described above specialises to a {\em homoclinic manifold}, for which the
results to be described also hold.  The intention is to characterise whether there are any persistent heteroclinic trajectories when
$ \eps \ne 0 $; that is, when the impulses are applied.  Of course, trajectories here are to be thought of in the sense described by
Lemma~\ref{lemma:diracunique}, in that all trajectories get reset when crossing $ t $ values in $ {\mathcal J} $.  It has 
already been established in Theorem~\ref{theorem:unstabledirac} that $ \Gamma $, which was originally a branch of the unstable
manifold of $ a $, with perturb to the unstable pseudo-manifold $ \tilde{\Gamma}_\eps^u(a) $.  Similarly, $ \Gamma $ when thought
of as a branch of the stable manifold of $ b $ will perturb by Theorem~\ref{theorem:stabledirac} to a stable pseudo-manifold
$ \tilde{\Gamma}_\eps^s(b) $.  Of course, there is no necessity for $ \tilde{\Gamma}_\eps^u(a) $ to coincide with 
$ \tilde{\Gamma}_\eps^s(b) $.  A picture of this situation in a time-slice $ t $ is shown in Figure~\ref{fig:distance}, where the
dashed curve is the unperturbed $ \Gamma $, with the unstable pseudo-manifold $ \tilde{\Gamma}_\eps^u(a) $ emanating from
$ a_+(t) $ and the stable pseudo-manifold  $ \tilde{\Gamma}_\eps^s(b) $ emanating from $ b_-(t) $.  The goal now is to express
the signed distance $ d(p,t,\eps) $, measured in the normal direction at $ \bar{x}(p) $ from $ x_\eps^s(p,t) $ to $ x_\eps^u(p,t) $,
in terms of the unperturbed flow and the spatial forms $ g_i $ associated with the impulses.

\begin{remark}[Standard Melnikov theory] \normalfont
The standard method for determining distances of this nature build on the Melnikov method \cite{melnikov,guckenheimerholmes,wiggins}, which in its original incarnation requires a steady two-dimensional area-preserving
flow possessing a $ \Gamma $ as in Figure~\ref{fig:distance}, to which is added a time-periodic perturbation.  However, both
area-preservation and time-periodicity can be relaxed \cite{tangential,siam_book}.  Thus, if the system were $ 
\dot{x} = f(x) + \eps g(x,t) $, 
with the $ \eps = 0 $ flow having identical hypothesis as in this article, but with $ g(x,t) $ being a bounded, sufficiently smooth
{\em function} as opposed to a distribution, then the Melnikov approach yields the fact that
\begin{equation}
\fl d(p,t,\eps) = \eps \frac{M(p,t)}{\left| f \left( \bar{x}(p) \right) \right|} + \O{\eps^2} \, , 
\label{eq:distance}
\end{equation}
where \cite{tangential,siam_book}
\begin{equation}
\fl M(p,t) =  \int_{-\infty}^\infty \! \! \exp \left[ \int_{\xi-t+p}^p {\mathrm{Tr}} \left[  D f \left( \bar{x}^u(\tau) \right) \right] \d \tau \right] f^\perp \left( 
\bar{x}^u(\xi \! - \! t \! + \! p) \right) \cdot g \left( \bar{x}^u(\xi\! -\! t \! + \!p), \xi \right) \, \d \xi \, .
\label{eq:melstandard}
\end{equation}
(Compare also with Remark~\ref{remark:impulseformal}, where the similar expression for only the unstable manifold is given.)
This simplifies to more familiar forms \cite{guckenheimerholmes,wiggins} under area-preserving flows in which
$ {\mathrm{Tr}} \, D f = 0 $.
\end{remark}

It is tempting to imagine that one can formally use (\ref{eq:melstandard}) when $ g(x,t) $ is a distribution, since the
integral is well-defined.  However, rigorously working through the integral equation shows that this is not quite the case:

\begin{theorem}[Distance between pseudo-manifolds]
\label{theorem:distance}
Let $ P $, $ T_u $ and $ -T_s $ be large, positive but finite, and suppose Hypotheses~\ref{hyp:fdirac}, \ref{hyp:gdirac} 
and \ref{hyp:heteroclinic} are satisfied.  Let $ p \in [-P,P] $ and $ t \in [T_s,T_u] \setminus {\mathcal J} $.  Then, the 
signed distance between $ \tilde{\Gamma}_\eps^u(a) $ and $ \tilde{\Gamma}_\eps^s(b) $ measured in the time-slice $ t $, at 
the location $ \bar{x}(p) $ in the direction $ f^\perp \left( \bar{x}(p) \right) $ is given by $ d(p,t,\eps) $ in (\ref{eq:distance}), 
where the Melnikov function is
\begin{equation}
\fl M(p,t) =  \sum_{i=1}^n  \left[ j_i(p,t) +  \int_{t_i}^t R_p(t-\xi) j_i(p,\xi) \, \d \xi \right] \, , 
\label{eq:melnikovdirac}
\end{equation}
where 
\begin{equation}
\fl j_i(p,t) = f^\perp \left( \bar{x}(t_i-t+p) \right) \cdot g_i \left(  \bar{x}(t_i-t+p) \right) \, , \, i=1,2, 3, \cdots, n \, ,
\label{eq:jump}
\end{equation}
and the resolvent $ R_p(t) $ is defined on $ \R \setminus \left\{ 0 \right\} $ by 
\begin{equation}
\fl R_p(t) = \left\{ \begin{array}{ll} {\mathcal L}^{-1} \left\{ \frac{\hat{F}_p^u(s)}{1 - \hat{F}_p^u(s)} \right\} (t) & ~~{\mathrm{if}}~ t > 0 \, , \\
{\mathcal L}^{-1} \left\{ \frac{\hat{F}_p^s(s)}{1 + \hat{F}_p^s(s)} \right\} (-t) &   ~~{\mathrm{if}}~ t < 0 \end{array} \right. \, , 
\label{eq:resolvent}
\end{equation}
with
\begin{equation}
\fl \hat{F}_p^u(s) = {\mathcal L} \left\{ {\mathrm{Tr}} \, D f \left( \bar{x}(p-t) \right) \right\}(s) \quad {\mathrm{and}} \quad 
\hat{F}_p^s(s) = {\mathcal L} \left\{ {\mathrm{Tr}} \, D f \left( \bar{x}(p+t) \right) \right\}(s) \, .
\label{eq:resolventtemp}
\end{equation}
\end{theorem}

\proof
See \ref{sec:distance}.
\proofend

\begin{theorem}[Heteroclinic persistence]
\label{theorem:heteroclinic}
Consider the conditions of Theorem~\ref{theorem:distance}.  If there exists $ (p_0,t_0) \in [-P,P] \times [T_s,T_u] \setminus
{\mathcal J} $ such that $ M(p_0,t_0) = 0 $ and $ \nabla M(p_0,t_0) \ne {\mathbf 0} $, then, for sufficiently small $ \left| \eps \right| $, 
there exists a $ (p,t) $ near $ (p_0,t_0) $ such that the trajectory of (\ref{eq:dirac}) passing through the time-slice $ t $ and
lying on the normal vector at $ \bar{x}(p) $ is heteroclinic: it approaches $ a $ in backwards time and $ b $ in forwards time.
\end{theorem}

\proof
This is a standard implicit function theorem argument which is no different from classical Melnikov results; see \cite[e.g.]{guckenheimerholmes}.
\proofend

\begin{remark}[Impulsive Melnikov function `is continuous'] \normalfont
\label{remark:melnikovdiraccontinuous}
Even though $ M(p,t) $ is defined for $ t \notin {\mathcal J} $, the expression (\ref{eq:melnikovdirac}) indicates
that for any $ t_i \in {\mathcal J} $, $ \lim_{t \uparrow t_i} M(p,t) = \lim_{t \downarrow t_i} M(p,t) $.  Thus the $ \left\{ t_i \right\} $
consist of removable singularities; if `filled in,' $ M $ would be continuous in $ t $.  The reason for this is that when crossing
a jump value $ t_i $, both $ x_\eps^u(p,t) $ and $ x_\eps^s(p,t) $ get reset according to the same jump map, which according to
Lemma~\ref{lemma:diracunique} is continuous.  Their relative distance to $ \O{\eps} $ turns out to be preserved during this
jump map; (\ref{eq:melnikovdirac}) implicitly establishes this fact.  This `continuity' of the Melnikov function in $ t $ was also
observed in an early attempt \cite{aperiodic} to rationalise flux under impulses (but
restricted to area-preservation).  However, if one obtains a zero for $ M $ at $ t $-values in $ {\mathcal J} $, this has no
physical interpretation in relation to Theorem~\ref{theorem:distance}.
\end{remark}

\begin{corollary}
\label{corollary:heteroclinic}
Under the conditions of Theorem~\ref{theorem:distance} suppose additionally that $ f $ is area-preserving.  Then the Melnikov
function (\ref{eq:melnikovdirac}) simplifies to
\begin{equation}
\fl M(p,t) =  \sum_{i=1}^n   j_i(p,t) \, , 
\label{eq:melnikovareapreserving}
\end{equation}
and moreover the conclusions of Theorem~\ref{theorem:heteroclinic} also hold.
\end{corollary}

\proof
Since $ {\mathrm{Tr}} D f = 0 $, the resolvent is zero, and the simplification is obvious.
\proofend

\begin{remark}[Formal Melnikov function for impulses] \normalfont
\label{remark:formalflux}
Remark~\ref{remark:impulseformal} has argued that the pseudo-manifold formul\ae{} for impulses
are equivalant to those obtained from the {\em smooth} Melnikov development by the formal substitution of Dirac
delta impulses into the relevant formul\ae{} (\ref{eq:melstandard}), {\em in the 
situation in which $ {\mathrm{Tr}} \, D f = 0 $}.  Since $ M(p,t) = M^u(p,t) - 
M^s(p,t) $, in area-preserving situations {\em only}, a formal Dirac delta substitution into the Melnikov function
(\ref{eq:melstandard}) does indeed yield the formula (\ref{eq:melnikovareapreserving}).  {\em This formal approach does not
work for non-area preserving flows}.
\end{remark}

\section{Example: heteroclinics in kicked Duffing oscillator}
\label{sec:duffing}

\begin{figure}[t]
\centering
\includegraphics[width = 0.6 \textwidth]{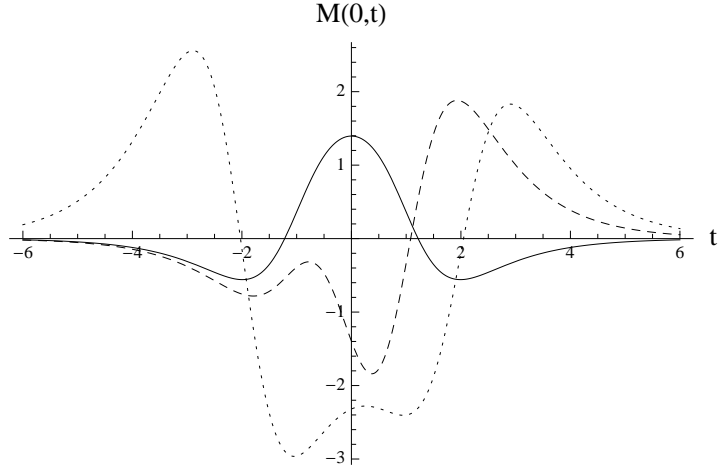}
\caption{The Melnikov function (\ref{eq:duffingmelnikov}) for the kicked Duffing oscillator for $ p = 0 $ for several choices of parameters:
$ n = 2 $, $ t_1 = -1 $, $ \gamma_1 = -1 $, $ t_2 = 1 $ and $ \gamma_2 = 1 $ [solid],  $ n = 3 $, $ t_1 = -1 $, $ \gamma_1 = -1 $, 
$ t_2 = 0 $, $ \gamma_2 = 1 $, $ t_3 = 1 $ and $ \gamma_3 = - 3 $ [dashed], and $ n = 2 $, $ t_1 = - 2 $, $ \gamma_1 = 3.7 $, $ t_2 = 2 $
and $ \gamma_2 = -2.7 $ [dotted].}
\label{fig:duffingmelnikov}
\end{figure}

Kicked oscillators are an oft-used paradigm in controlling chaos \cite{linyoungkick,linyoungshear,franaszek,bapat}, and here the Duffing oscillator \cite{franaszek,bapat,holmes,holmeswhitley,wiggins,wigginsduffing,mancho,jimenez,tangential,belleyvirgilio_duffing} is chosen. If subject to a finite number of kicks at times $ \left\{ t_1, t_2, \cdots , t_n \right\} $, the undamped impulsively-forced Duffing oscillator is
given by
\begin{equation}
\fl \ddot{x} - x + x^3 = \eps \sum_{i=1}^n \gamma_i \delta(t - t_i)
\label{eq:duffing}
\end{equation}
where the $ \gamma_i \in \R $ represent the sizes and directions of the kicks, $ 0 < \eps \ll 1 $, $ x \in \R $, and the overdot
represents the time-derivative.  In this case there is no ambiguity in writing the evolution in terms of a differential equation
since the impulsive terms are spatially-independent; the fact that the vector field is area-preserving (as will be seen) renders
this approach particularly attractive.   To be consistent with the notation of this article, set $ x = x_1 $ and $ \dot{x} = x_2 $, to
get the system
\begin{equation}
\fl \frac{d}{d t} \left( \begin{array}{c} x_1 \\ x_2 \end{array} \right) =
\left( \begin{array}{c} x_2 \\ x_1 - x_1^3 + \eps \sum_{i=1}^n \gamma_i \delta(t - t_i) \end{array} \right) \, .
\label{eq:duffingsystem}
\end{equation}
When  $ \eps = 0 $, the phase portrait of the Duffing oscillator is well-known to have a figure-eight structure in the
$ x_1 x_2 $-plane centred at the saddle point at the origin, with the two rings of the figure-eight each representing a heteroclinic
connection \cite[e.g.]{holmes,wiggins,tangential}.  The right branch is representible as a solution to (\ref{eq:duffing}) with $ \eps = 0 $ by
\begin{equation}
\fl \left( \begin{array}{c} \bar{x}_1(t) \\ \bar{x}_2(t) \end{array} \right) = \left( \begin{array}{c} \sqrt{2} \sech t \\
- \sqrt{2} \sech t \tanh t \end{array} \right) \quad , \quad t \in \R \, .
\label{eq:duffingheteroclinic}
\end{equation}
Now, in this case 
\[
\fl f^\perp \left( \bar{x}(t) \right) = \left( \begin{array}{c} \bar{x}_1(t)^3 - \bar{x}_1(t) \\ \bar{x}_2(t) \end{array} \right) = 
\left( \begin{array}{c} \sqrt{2} \sech t \left[ 2 \sech^2 t - 1 \right] \\ \sqrt{2} \sech t \tanh t \end{array} \right) \, , 
\]
and so from (\ref{eq:jump}), 
\[
\fl j_i(p,t) =  \gamma_i \sqrt{2} \sech \left( t_i - t +p\right) \tanh \left( t_i - t +p\right) \, .
\]
Now in this case $ {\mathrm{Tr}} \, D f = 0 $, and so Corollary~\ref{corollary:heteroclinic}
can be used directly.  The Melnikov function is therefore
\begin{equation}
\fl M(p,t) = \sqrt{2} \sum_{i=1}^n  \gamma_i \sech \left( t_i - t +p\right) \tanh \left( t_i - t +p\right) \, .
\label{eq:duffingmelnikov}
\end{equation}
With the choice $ p = 0 $, the distance between the perturbed pseudo-manifolds will be measured at $ (x_1,x_2) = (\sqrt{2},0) $,
along the normal direction $ (1,0) $.  If $ M(0,t) $ has a simple zero at $ t $, then in the time-slice
$ t $ there will be a heteroclinic trajectory passing near the point $ (\sqrt{2},0) $.  Now with $ p = 0 $, each term in 
(\ref{eq:duffingmelnikov}) is odd about $ t_i $,
so for example if $ n = 1 $, the presence of a simple zero at $ t_1 $ can be immediately imputed.  However, this is in the set
$ {\mathcal J} $, and therefore one {\em cannot} automatically 
conclude the presence of persistent heteroclinics if $ n = 1 $ and $ \gamma_1 \ne 0 $.  With this in mind, Figure~\ref{fig:duffingmelnikov} shows
the function (\ref{eq:duffingmelnikov}) for several choices of  $ n $, $ t_i $ and $ \gamma_i $, with the zeros of $ M(0,t) $ 
in each instance indicating the presence of heteroclinic trajectories which backwards and forwards asymptote to the rest
state $ (x_1,x_2) = (0,0) $.  The zeros visible in all cases are simple and removed from $ {\mathcal J} $.

\section{Transport and flux}
\label{sec:flux}

\begin{figure}[t]
\centering
\includegraphics[width=0.72\textwidth]{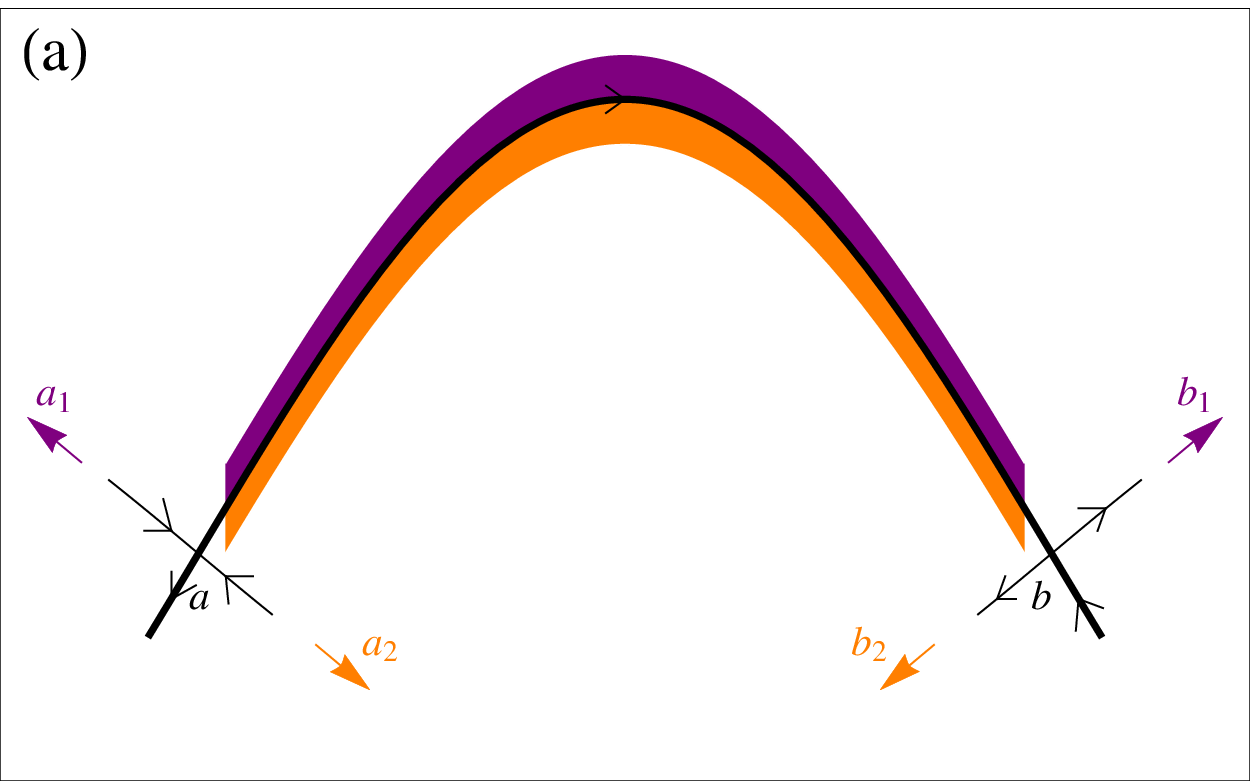} \\
\includegraphics[width=0.72\textwidth]{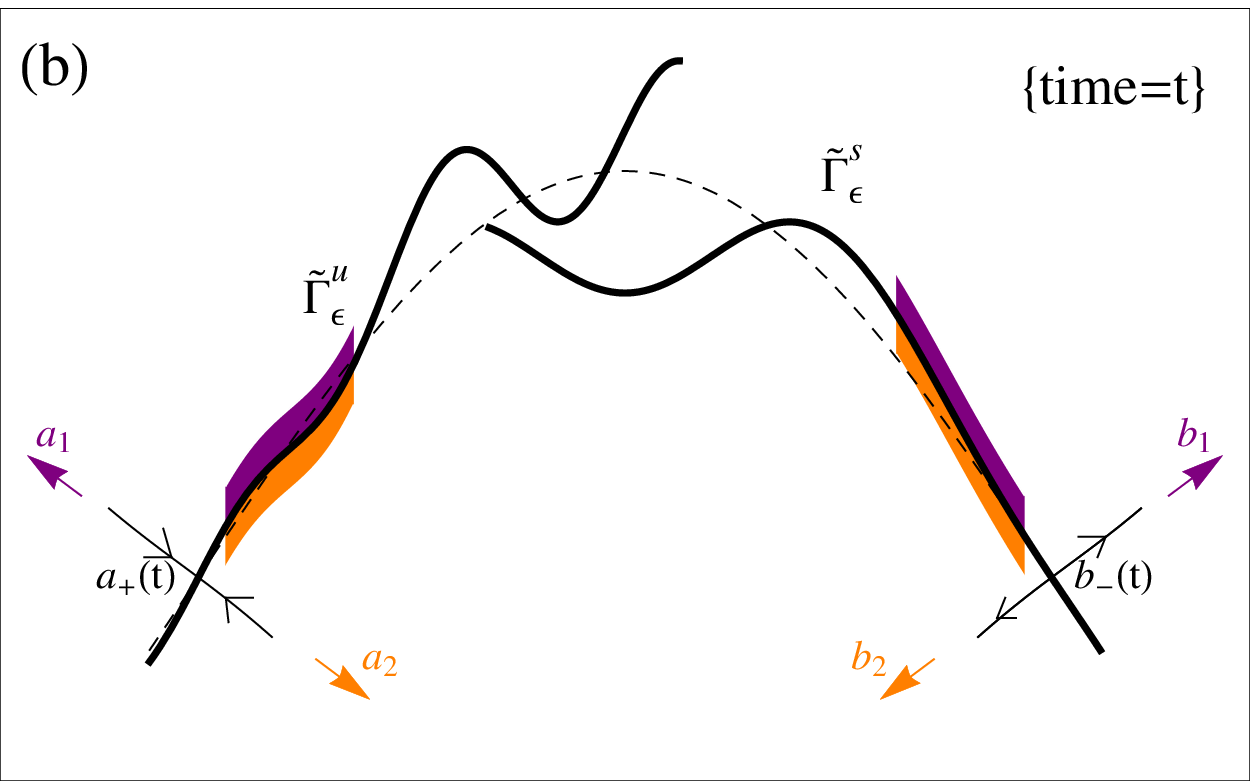} \\
\includegraphics[width=0.72\textwidth]{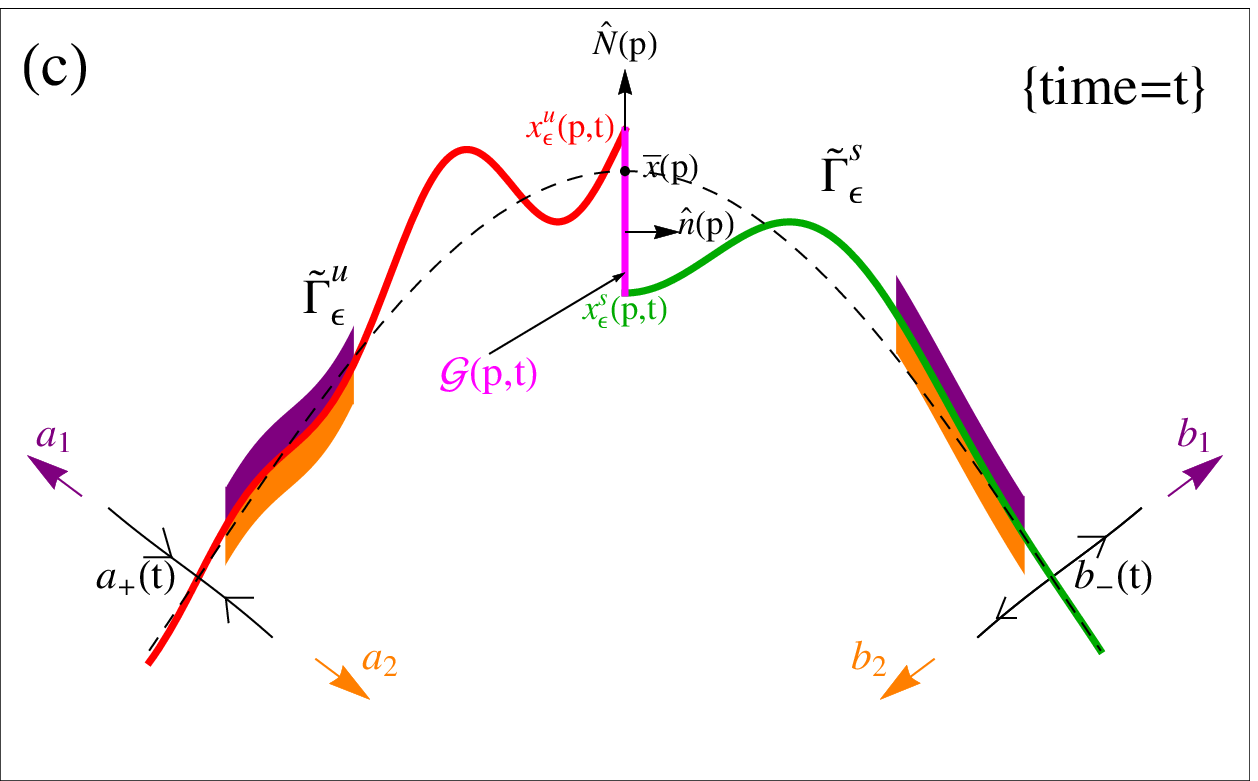}
\caption{The flow-separating pseudo-manifolds when (a) $ \eps = 0 $ and (b) $ \eps \ne 0 $, with the pseudo-separatrix
construction of Definition~\ref{definition:pseudoseparatrix} shown in (c).}
\label{fig:flux}
\end{figure}

An unbroken codimension-$1$ heteroclinic manifold is an important flow barrier in autonomous flows; trajectories on the opposite sides experience different fates.
This is easily seen by considering Figure~\ref{fig:flux}(a), which shows the $ \eps = 0 $ heteroclinic manifold along with
a shaded strip of nearby `particles' lying on both sides of the manifold.  In forwards time, the upper (darker) collection will
get pulled away in the direction indicated by the vector $ b_1 $, which is associated with one branch of the unstable manifold of
$ b $.  These particles will be termed `$b_1 $-forward' particles. On the other hand, the lighter group, lying below the heteroclinic manifold, will get pulled away from $ b $ in the {\em opposite} direction indicated by $ b_2 $, representing the opposite branch of $ b $'s unstable manifold.  These are `$b_2 $-forward' particles.  If now considering the
fate of each of these two groups in backwards time, the upper (darker) group will get pulled away from $ a $ in the direction
$ a_1 $ (`$a_1 $-backward' particles), while the lower (lighter) group will experience repulsion from $ a $ in the direction $ a_2 $
(`$a_2 $-backward' particles).   What is clear in this instance is that the $ a_1 $-backward particles are identical to the
$ b_1 $-forward particles, and also the $ a_2 $-backward particles are the same as the $ b_2 $-forward ones. The clear distinction between these groups, which are divided by the heteroclinic manifold, highlights the idea
that the heteroclinic manifold is a flow separator.  In this instance, also note that there is no flux of particles from
one group to the other; the manifold is impermeable. 

When $ \eps \ne 0 $, the two pseudo-manifolds $ \tilde{\Gamma}_\eps^u(a) $ and $ \tilde{\Gamma}_\eps^s(b) $ are created,
which do not need to coincide as is shown in Figure~\ref{fig:distance}.  Moreover, these evolve with time. 
A possible situation, in a time-slice $ t \in [T_s,T_u] \setminus {\mathcal J} $, is shown
in Figure~\ref{fig:flux}(b), where the dashed curve is the heteroclinic manifold $ \Gamma $, as shown in Figure~\ref{fig:flux}(a).  To avoid clutter, the pictured unstable and stable pseudo-manifolds, which emanate respectively from $ a_+(t) $
and $ b_-(t) $, have been clipped after proceeding some distance away from these points.  However, these pseudo-manifolds
could be very complicated, intersecting each other in various ways.  Now, the coloured particle
groups on the two sides of $ \tilde{\Gamma}_\eps^u $,  will in backwards time always remain on the two sides
of the time-evolving $ \tilde{\Gamma}_\eps^u $.  In going back in time, each time the (backward) jump map operates, 
the pseudo-manifold will get reset, {\em and preserve the fact that these particle groups are on its two sides}.  Once gone back in time beyond $ t_1 $---the first impulse time---these particles will
therefore remain separated by the original unperturbed unstable manifold $ \Gamma $.  Thus, the darker group of particles will
be $ a_1 $-backward, whereas the lighter group will be $ a_2 $-backward.  The particles on
the two sides of $ \tilde{\Gamma}_\eps^u $ will therefore in backwards time get pushed apart in the directions $ a_1 $ and $ a_2 $,
at an exponential rate.  So, even when $\eps \ne 0 $, 
$ \tilde{\Gamma}_\eps^u $ continues to be a flow separator in backward time.  The same argument follows for the groups of particles shown
near $ b_-(t) $ on the two sides of $ \tilde{\Gamma}_\eps^s $ in forward time: the darker group will be $ b_1 $-forward,
the lighter will be $ b_2 $-forward, and these will be pushed apart exponentially.

The difference between Figure~\ref{fig:flux}(a) and \ref{fig:flux}(b) is that the separation between the lighter and darker
groups of particles no longer occurs in {\em both} forward and backward time.  So for the situation pictured in Figure~\ref{fig:flux}(b), some of the lighter $ a_2 $-backward particles near $ a_+(t) $ may in {\em forward} time be $ b_1 $-forward.  Whether this happens or not depends on the
particular intersection pattern between $ \tilde{\Gamma}_\eps^u $ and $ \tilde{\Gamma}_\eps^s $, enabling a {\em transfer} from
one side to the other.  These pseudo-manifolds may intersect in
complicated ways, or not at all, and depending on this, some particles will be $ b_1 $-forward, while others are
$ b_2 $-forward.  It is because of this that one can imagine that {\em transport} has occurred across the
barrier $ \Gamma $, which was impermeable in backward and forward time when $ \eps = 0 $.  The difficulty now is in {\em quantifying} the transport occurring as a result of the broken heteroclinic, bearing in mind the possibility that $ \tilde{\Gamma}_\eps^u $ and $ \tilde{\Gamma}_\eps^s $ may intersect (or not) in numerous ways.

This issue is not confined to impulsive perturbations: the same problem arises 
even with smooth, aperiodic, perturbations which result in smooth stable and unstable manifolds $ \Gamma_\eps^u(a) $ and $ \Gamma_\eps^s(b) $.  The resolution to this is to consider a {\em time-varying flux} of particles between the lighter and
the darker
groups of particles \cite{aperiodic,open}.  The description for this smooth situation is not any different from the current
situation, except that in this impulsive case, the fact that there is no perturbation for $ t < t_1 $ or $ t > t_n $ makes thing
simpler in some senses, since $ \tilde{\Gamma}_\eps^u = \Gamma $ for $ t < t_1 $, and $ \tilde{\Gamma}_\eps^s = \Gamma $
for $ t > t_n $.   Consider Figure~\ref{fig:flux}(b), which shows
that the perturbed stable and unstable pseudo-manifolds do not coincide.  To make sense of a fluid transfer between the
two groups of particles, fix a $ p \in [-P,P] $, and consider a point $ \bar{x}(p) $ 
on the unperturbed heteroclinic manifold.  Draw a perpendicular vector to $ \Gamma $ at this point, i.e., in the direction defined
by $ \hat{N}(p) $.  This will intersect
$ \tilde{\Gamma}_\eps^u $ at $ x_\eps^u(p,t) $, and $ \tilde{\Gamma}_\eps^s $ at $ x_\eps^s(p,t) $, which are both 
in a $ {\mathcal O}(\eps) $-neighborhood, since the pseudo-manifolds
have perturbed within such a distance.  This construction is shown in Figure~\ref{fig:flux}(c), based on the geometry
of Figure~\ref{fig:flux}(b).  

\begin{definition}[Pseudo-separatrix]
\label{definition:pseudoseparatrix}
The {\em pseudo-separatrix} $ {\mathcal Q}(p,t) $ is the union of three curves, as shown in Figure~\ref{fig:flux}(c):
\begin{enumerate}
\item The unstable pseudo-manifold curve $ \tilde{\Gamma}_\eps^u $ [red] emanating from $ a_+(t) $ until it reaches 
$ x_\eps^u(p,t) $, where the normal vector to $ \Gamma $ drawn at $ \bar{x}(p) $ intersects it;
\item The stable pseudo-manifold curve  $ \tilde{\Gamma}_\eps^s $ [green] emanating from $ b_-(t) $ until it reaches 
$ x_\eps^s(p,t) $, where the normal vector to $ \Gamma $ drawn at $ \bar{x}(p) $ intersects it;
\item The line [magenta] which connects these two curves along the normal vector at $ \bar{x}(p) $, which shall
be called the {\em gate} $ {\mathcal G}(p,t) $.
\end{enumerate}
\end{definition}

The curve $ {\mathcal Q}(p,t) $ is of course not a pure flow separator when $ \eps \ne 0 $.  Think now of its evolution with $ t $, bearing in mind
that the unstable/stable pseudo-manifold segments will be evolving with time, and the gate will have to be extended/shrunk depending on the locations of $ x_\eps^{s,u}(p,t) $.  
Viewing Figures~\ref{fig:flux}(b) and \ref{fig:flux}(c)
together, the following observations can be made:
\begin{itemize}
\renewcommand{\labelitemi}{$-$}
\item Particles below $ {\mathcal Q} $ will be $ a_2 $-backward  since below $ \tilde{\Gamma}_\eps^u(a) $;
\item Particles above $ {\mathcal Q} $ will be $ b_1 $-forward time since above $ \tilde{\Gamma}_\eps^s(b) $;
\item Over an infinitesimal time, no particles will cross either the $ \tilde{\Gamma}_\eps^u(a) $ or the  $ \tilde{\Gamma}_\eps^s(b) $
segments of $ {\mathcal Q} $ since they are material curves evolving with time;
\item Thus, the only transfer from $ a_2 $-backward to $ b_1 $-forward particles can occur via particles instantaneously flowing through the gate.
\end{itemize}

The task now is to define the {\em instantaneous flux} of particles across $ {\mathcal Q}(p,t) $ from $ a_2 $-backward to $ b_1 $-forward.  Let this be denoted
by $ \phi(p,t) $, with $ p $ denoting the location of the gate, and $ t \in [T_s,T_u] \setminus {\mathcal J} $, time.  By the above
argument, $ \phi(p,t) $ is therefore the instantaneous flux across just the gate.  To express this, let $ \ell $ be an arclength
parametrisation of the line segment $ {\mathcal G}(p,t) $, chosen such that $ \ell = 0 $ at 
$ x_\eps^s(p,t) $, and $ \ell = L(p,t) > 0 $ at the other endpoint
$ x_\eps^u(p,t) $.  Use $ f(\ell) $ as the short-hand notation for the instantaneous velocity at a location $ \ell $ on
the gate (since $ t \notin {\mathcal J} $ the velocity only contains the unperturbed component $ f $).  Let $ \hat{n}(\ell) $ be 
the unit normal vector to the gate at a general location, chosen with direction consistent with $ f (\ell) $ (i.e., 
consonant with $ f \left( x_\eps^{u,s}(p,t) \right) \cdot \hat{n} > 0 $).  Then, the instantaneous
flux is defined by
\begin{equation}
\fl \phi(p,t) := \int_0^{L(p,t)}
f \left( \ell \right) \cdot \hat{n}(\ell) \, \d \ell \, .
\label{eq:fluxdef}
\end{equation}
As defined, this gives precisely a quantity of fluid per unit time, crossing $ {\mathcal Q} $, and the flux
depends on the gate location (parametrised by $ p $) and time $ t $.
Now, in the situation pictured in Figure~\ref{fig:flux}(c), one would get a positive instantaneous flux, which is 
therefore associated with a transfer from $ a_2 $-backward to $ b_1 $-forward particles.  The impact of this, in relation
to transport across $ \Gamma $, is that the transport occurs from the lower to the upper fluids, corresponding to a direction
$ \hat{N}(p) $.  Thus, a positive $ \phi $ implies instantaneous transport
across $ \Gamma $ in the direction of $ + \hat{N}(p) $, whereas a negative $ \phi $ is associated with transport in the
direction of $ - \hat{N}(p) $.  The latter case occurs if  the stable pseudo-manifold met $ {\mathcal G} $
at a higher point than does the unstable pseudo-manifold, and then the transfer is from the upper to the lower fluid instead, 
i.e., from $ a_1 $-backward to $ b_2 $-forward.  At instances in which $ \tilde{\Gamma}_\eps^u $
and $ \tilde{\Gamma}_\eps^s $ intersect {\em exactly} on $ {\mathcal G}(p,t) $, the instantaneous flux is zero.

\begin{theorem}[Flux]
The instantaneous flux across the pseudo-separatrix $ {\mathcal Q}(p,t) $, associated with a gate location $ \bar{x}(p) $ and a time $ t \in [T_s,T_u] \setminus {\mathcal J} $, is 
given by
\begin{equation}
\fl \phi(p,t) = \eps M(p,t) + \O{\eps^2} \, .
\label{eq:flux}
\end{equation}
\end{theorem}

\proof
The velocity at all points on $ {\mathcal G} $ is given by $ f (\ell) = f \left( \bar{x}(p) \right) + {\mathcal O}(\eps) $, since
all points on $ {\mathcal G} $ are $ {\mathcal O}(\eps) $-close to $ \bar{x}(p) $.    The normal vector $ \hat{n}(\ell) $ 
is also to leading-order equal to the unit normal in the direction of $ f \left( \bar{x}(p) \right) $.  However, Theorem~\ref{theorem:distance} establishes that the
leading-order displacement 
\[
\fl x_\eps^u(p,t) - x_\eps^s(p,t) = \eps \frac{M(p,t)}{\left| f \left( \bar{x}(p) \right) \right|} + {\mathcal O}(\eps^2) \, , 
\]
and thus
\begin{eqnarray*}
\fl \phi(p,t) & = & \int_0^{L(p,t)}
\left[ f \left( \bar{x}(p) \right) + {\mathcal O}(\eps) \right] \cdot \left[ \frac{f \left( \bar{x}(p) \right)}{\left|
f \left( \bar{x}(p) \right) \right|} + {\mathcal O}(\eps) \right] \, \d \ell \\
& = & \left| f \left( \bar{x}(p) \right) \right| \int_0^{L(p,t)} \d \ell + {\mathcal O}(\eps^2) \\
& = &  \left| f \left( \bar{x}(p) \right) \right| \left[ \eps \frac{M(p,t)}{\left| f \left( \bar{x}(p) \right) \right|} + {\mathcal O}(\eps^2)
\right] + {\mathcal O}(\eps^2) \\
& = & \eps M(p,t) + {\mathcal O}(\eps^2) \, , 
\end{eqnarray*}
where the fact that the length of the gate is $ {\mathcal O}(\eps) $ has been used at the second step.
\proofend

Thus, the leading-order instantaneous flux, as a time-varying entity, {\em is} the Melnikov function (this result, valid for
general smooth time-varying perturbations \cite{aperiodic,open}, generalises thinking of the integral of the Melnikov function 
as a measure of  lobe-dynamics transport in time-periodic flows \cite{romkedar}). Basically, the general time development for smooth
perturbations \cite{aperiodic,open} applies to this impulsive setting as well, with the understanding that 
a positive flux at some time $ t $ implying instantaneous 
transfer across $ \Gamma $ in the direction of $ + \hat{N}(p) $.

\section{Example: flux in an eddy due to an underwater explosion}
\label{sec:eddy}

\begin{figure}[t]
\centering
\includegraphics[width=0.45 \textwidth]{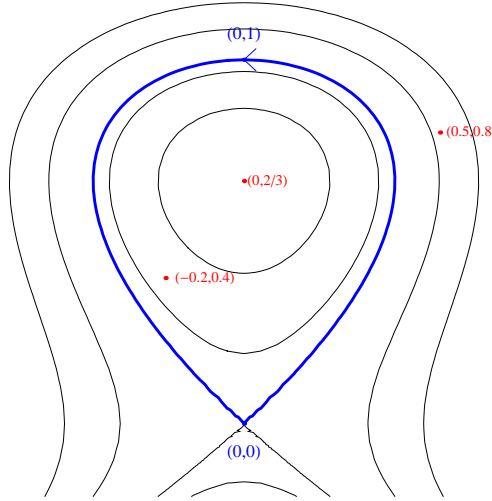}
\caption{The eddy structure of (\ref{eq:eddy}), in which the effect of an underwater explosion centred at each of the red
dots will be assessed.}
\label{fig:eddy}
\end{figure}

The simplest nontrivial situation in which the results of the previous section are applicable will be first considered,
with a more complex example provided in the subsequent section.  Thus, the flow will be area-preserving, and there shall be
only one time at which an impulse applies.  This shall be in the context of a highly idealised situation in which the flux in an 
oceanic eddy due to an underwater explosion is to be assessed.    In the absence of an explosion, the
model for the oceanic eddy shall be given by
\begin{equation}
\fl \frac{d}{d t} \left( \begin{array}{c} x_1 \\ x_2 \end{array} \right) = 
\left( \begin{array}{c} 2 x_2 - 3 x_2^2\\ 2 x_1  \end{array} \right)  \, , 
\label{eq:eddy}
\end{equation}
whose phase portrait is shown in Figure~\ref{fig:eddy}.   This is a kinematic model for an oceanic eddy; in this form this
models for example a warm-core eddy detaching northward from the Gulf Stream \cite{viscous,eddy,delcastillonegretemorrison}.
The `outermost' closed loop, shown in blue, is a homoclinic
trajectory $ \Gamma $ associated with the point $ (0,0) $, and can be represented by $ x_1 = \pm x_2 \sqrt{
1 - x_2} $, $ x_2 \in (0, 1] $.  It is across this that the flux due to an underwater explosion centred at
an arbitrary point  $ (\tilde{x}_1,\tilde{x}_2) $, either inside or outside the eddy, is to be assessed.  
Using a symmetric time-parametrisation ensuring that $ t = 0 $ corresponds to
the top-most point $ (0,1) $, the homoclinic trajectory can be obtained as
\[
\fl \bar{x}(t) = \left( \begin{array}{c} \bar{x}_1(t) \\ \bar{x}_2(t) \end{array} \right) = 
\left( \begin{array}{c} -\sech^2 t \tanh t \\ 
\sech^2 t \end{array} \right) \,  ,
\]
and consequently,
\[
\fl f^\perp \left( \bar{x}(t) \right) =  \left( \begin{array}{c} 
2 \sech^2 t \tanh t \\ 
\left[ \cosh \left( 2 t \right) - 2 \right] \sech^4 t \end{array} \right) \, .
\]
Now, an underwater explosion is assumed to occur at time $ t_1 = 0 $, at the location $ (\tilde{x}_1,\tilde{x}_2) $.   
This can be anywhere in the fluid, but not on the homoclinic, and thus $ \tilde{x}_1  \ne \pm \tilde{x}_2 \sqrt{1 - \tilde{x}_2} $ if $ \tilde{x}_2 \in [0,1] $.  A plausible model is that the explosion generates an impulsive velocity
radially outwards from $ (\tilde{x}_1,\tilde{x}_2) $, and that the effect of this diminishes with the distance from this point.  Thus, suppose that
\[
\fl g_1 \left( x \right) = \frac{1}{\left( x_1 - \tilde{x}_1 \right)^2 + \left( x_2 - \tilde{x}_2 \right)^2} \left( \begin{array}{c} 
x_1 - \tilde{x}_1 \\ x_2 - \tilde{x}_2 \end{array} \right) \, ,
\]
bearing in mind that the resulting velocity is this multiplied by $ \eps \delta(t) $.

\begin{figure}[t]
\centering
\includegraphics[width=0.6 \textwidth]{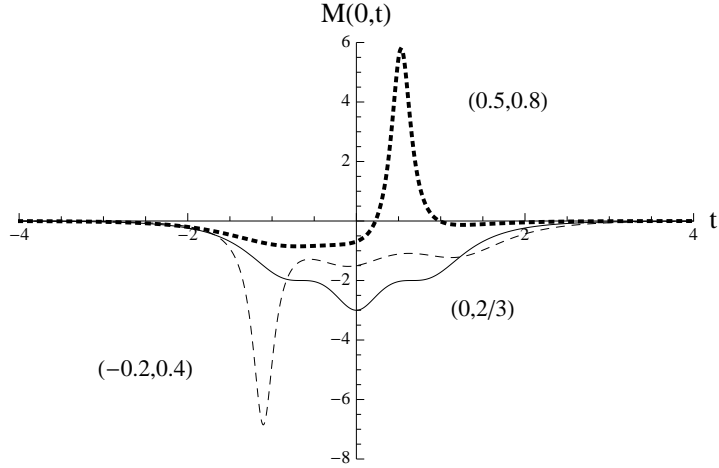}
\caption{Leading-order flux functions for the eddy of Figure~\ref{fig:eddy} associated with explosions centred at each
of the red dots in Figure~\ref{fig:eddy}, computed according to (\ref{eq:melexplosion}).}
\label{fig:melexplosion}
\end{figure}

Choose $ p = 0 $; the gate is therefore located at the uppermost point $ (0,1) $ on the homoclinic.  Using Corollary~\ref{corollary:heteroclinic}, the resulting flux is therefore $ \phi(0,t) = \eps M(0,t) + {\mathcal O}(\eps^2) $, where $ M(0,t) = f^\perp
\left( \bar{x}(-t) \right) \cdot g_1 \left( \bar{x}(-t) \right) $, and hence
\begin{equation}
\fl M(0,t) = \frac{- 2 \sech^2 t \tanh t \left( \sech^2 t \tanh t \! - \! \tilde{x}_1 \right) + \left[ \cosh \left( 2 t \right) \! - \! 2 \right]
\sech^4 t \left( \sech^2 t \! - \! \tilde{x}_2 \right)}
{ \left( \sech^2 t \tanh t - \tilde{x}_1 \right)^2 + \left( \sech^2 t - \tilde{x}_2 \right)^2} \, .
\label{eq:melexplosion}
\end{equation}
A positive $ M $ would indicate flux into the eddy, with negative $ M $ flux out of it.  Thus, $ M $ can in this instance be
thought to represent precisely the rate of change of the size of the eddy.

The function $ M $ computed for each
of the three locations of the explosion shown by a red dot in Figure~\ref{fig:eddy} are shown in Figure~\ref{fig:melexplosion},
with the explosion locations stated in the form $ (\tilde{x}_1,\tilde{x}_2) $ adjacent to each curve.  The explosion occurring at the 
centre $ (0,2/3) $ of the eddy results in fluid leaving the eddy at all times, with the flux decaying as $ t \rightarrow \pm \infty $ (as it must in
all cases).  The total fluid leaving the eddy can be computed by the area between the curve and the $ t $-axis in Figure~\ref{fig:melexplosion}, i.e., the integral of (\ref{eq:melexplosion}) over $ \R $.  Therefore, the result of this particular explosion on the eddy
is that it diminishes in size by losing its warm interior waters to the outer colder sea.  The explosion occurring at the exterior location $ (0.5,0.8) $ results first in fluid leaving the eddy, and then later in a pulse of fluid entering the eddy at around $ t \approx 0.5 $. 

\begin{figure}
\begin{center}
\includegraphics[width=0.4 \linewidth, height=0.3 \textheight]{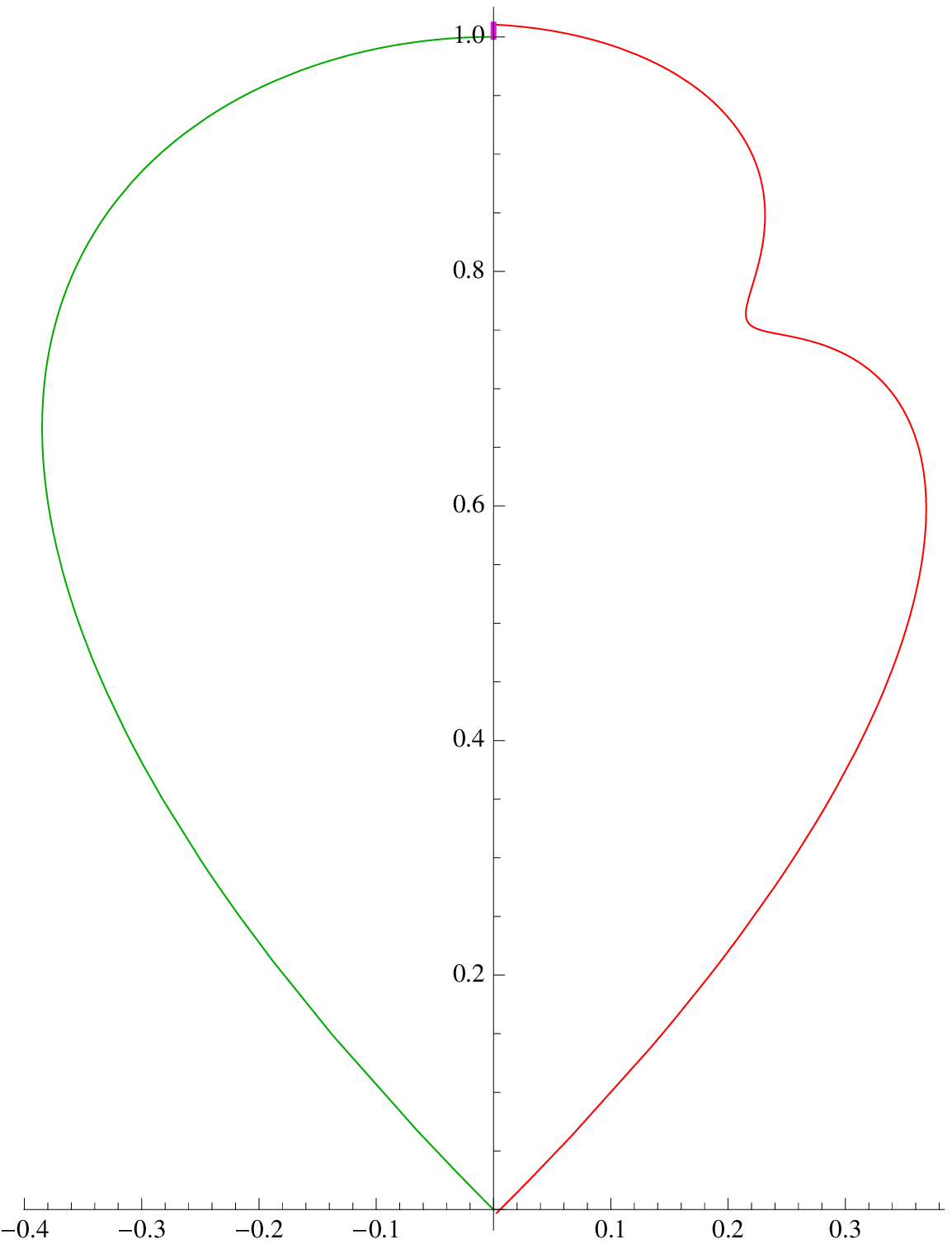} \hfill
\includegraphics[width=0.4 \linewidth, height=0.3 \textheight]{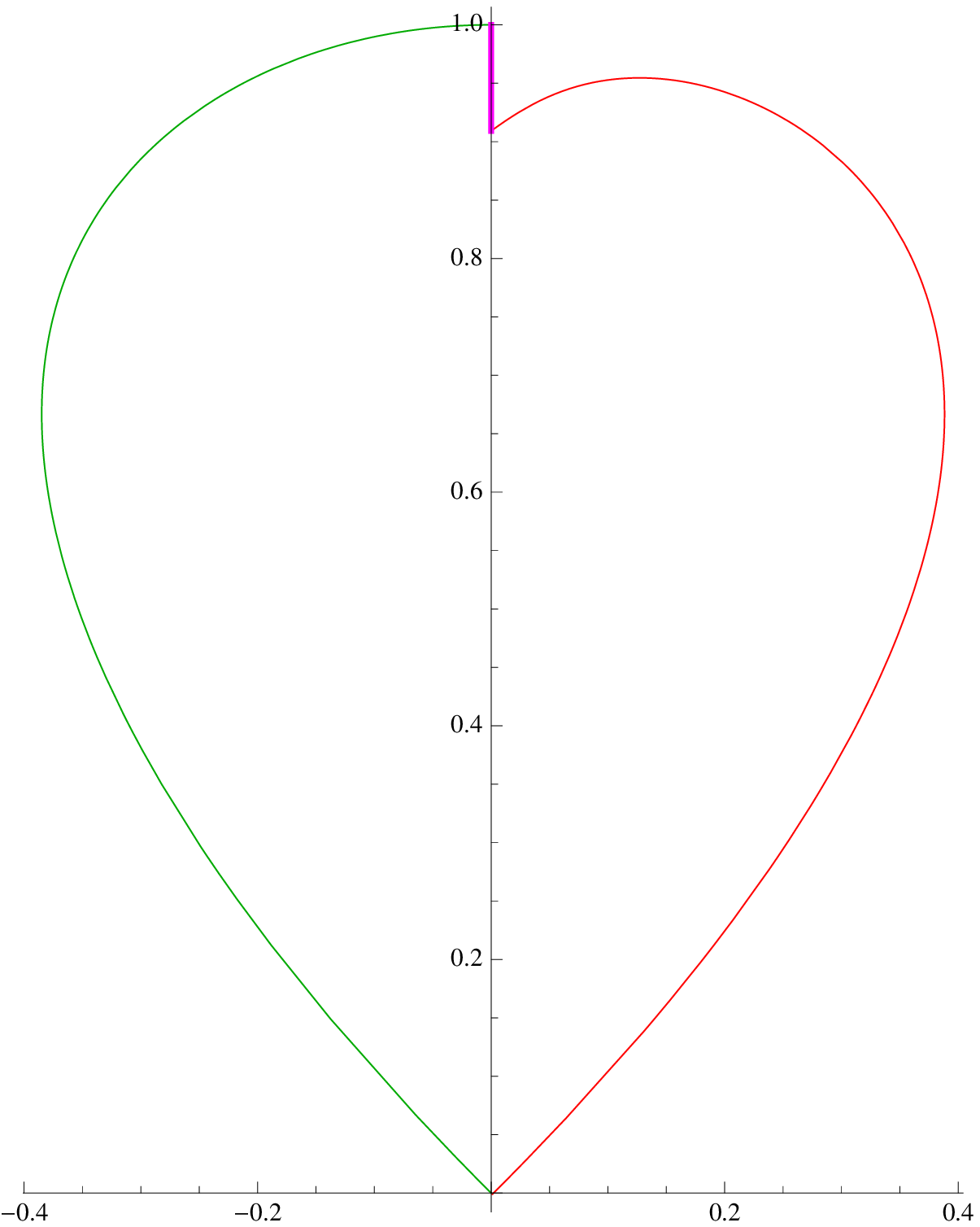} 
\caption{The pseudo-separatrix formed by the stable pseudo-manifold [green], the unstable pseudo-manifold [red] and the
gate [magenta] for the explosion centred at $ (0.5,0.8) $ and $ \eps = 0.02 $ using (\ref{eq:manifoldexplosion}) at times $ t = 0.1 $ [left] and
$ t = 0.6 $ [right].}
\label{fig:pseudoexplode}
\end{center}
\end{figure}

For this particular example, approximations for the pseudo-manifolds, and the pseudo-separatrix $ {\mathcal Q} $, can be
explicitly constructed using Corollaries~\ref{corollary:unstabledirac} and \ref{corollary:stabledirac}.  The pseudo-separatrix
is in fact a nominal boundary to the eddy, in an instance in which an absolute boundary does not exist.  Moreover, it is chosen
in such a way as to enable the quantification of waters into or out of the eddy, using the Melnikov function as this leading-order
flux.  Using Theorem~\ref{theorem:unstabledirac}, the
unstable pseudo-manifold emanating from $ (0,0) $ to the gate (at $ p = 0 $), in the time-slice $ t $ would be given in the form
\begin{equation}
\fl x_\eps^u(p,t) \approx \left( \begin{array}{c} \! \! \! - \sech^2 p \tanh p \\ \sech^2 p \end{array} \! \! \right) \! + \!  
\frac{\eps M^u(p,t)}{4 \sech^4 p \tanh^2 p \! + \! \left[ \cosh \left( 2 p \right) \! - \! 2 \right]^2 \! \sech^8 p} 
\left( \begin{array}{c} \! \! \! 
2 \sech^2 p \tanh p \\ 
\left[ \cosh \left( 2 p \right) \! - \! 2 \right] \! \sech^4 p \end{array} \! \! \right) 
\label{eq:manifoldexplosion}
\end{equation}
where
\[
\fl M^u(p,t) = \I_{(0,\infty)}(t) \Lambda(p-t) \quad {\mathrm{for}} \quad p < 0 \, , 
\]
in which
\[
\fl \Lambda(\xi) := \frac{2 \sech^2 \xi \tanh \xi \left( - \sech^2 \xi \tanh \xi - \tilde{x}_1 \right) + \left[ \cosh \left( 2 \xi \right) - 2 \right]
\sech^4 \xi \left( \sech^2 \xi - \tilde{x}_2 \right)}
{ \left( - \sech^2 \xi \tanh \xi - \tilde{x}_1 \right)^2 + \left( \sech^2 \xi - \tilde{x}_2 \right)^2} \, .
\]
The restriction $ p < 0 $ ensures that this pseudo-manifold is only drawn from the point $ (0,0) $ until it intersects the gate drawn
at $ (0,1) $.
The stable pseudo-manifold would be given by the expression (\ref{eq:manifoldexplosion}) with the superscript $ u $ replaced
by $ s $, and where
\[
\fl M^s(p,t) = - \I_{(-\infty,0)}(t) \Lambda(p-t) \quad {\mathrm{for}} \quad p > 0 \, . 
\]
The gate would connect the stable and unstable pseudo-manifolds, and together these would form the pseudo-separatrix
across which the flux is assessed; (\ref{eq:melexplosion}) is the leading-order expression for this.  The pseudo-separatrices
formed by these expressions are plotted for the explosion centred at $ (0.5,0.8) $, at two different times, in Figure~\ref{fig:pseudoexplode}.  The colour-coding red/green/magenta associated with the three curves comprising the
pseudo-seperatrix as outlined in Definition~\ref{definition:pseudoseparatrix} has been followed.  Notice the impact of the explosion at $ t = 0 $ has had a dramatic impact on the
unstable pseudo-manifold at $ t = 0.1 $, pushing it away from the explosion centre $ (0.5,0.8) $.
The relative positioning of the
stable and unstable pseudo-manifolds indicates that the flux is negative (out of the eddy)
at $ t = 0.1 $, because the unstable pseudo-manifold is slightly higher than the stable one on the gate.  Similarly, the right panel of Figure~\ref{fig:pseudoexplode} shows that the flux 
is positive at $ t = 0.6 $, with a significantly larger magnitude because of the larger gate.  Thus the eddy is instantaneously shrinking (slightly) at $ t = 0.1 $, but expanding at $ t = 0.6 $.  These observations are consistent with the dotted curve in
Figure~\ref{fig:melexplosion}, which shows the flux variation for this situation.  
  A small value of $ \eps $ was
needed in producing these plots because the function $ g $ itself has a singularity at the explosion centre, and thus is very large
on the homoclinic if the centre is near to it.

\section{Example: flux in an expanding flow}
\label{sec:expanding}

\begin{figure}[t]
\centering
\includegraphics[width=0.6 \textwidth, height=0.3\textheight]{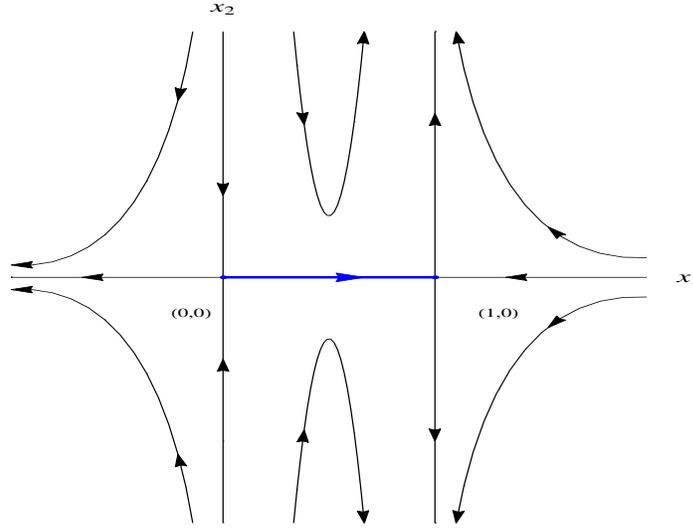}
\caption{Phase plane associated with (\ref{eq:madeup}), for the flux computation example.}
\label{fig:madeupphase}
\end{figure}

\begin{figure}[t]
\includegraphics[width=0.47 \textwidth, height=0.3\textheight]{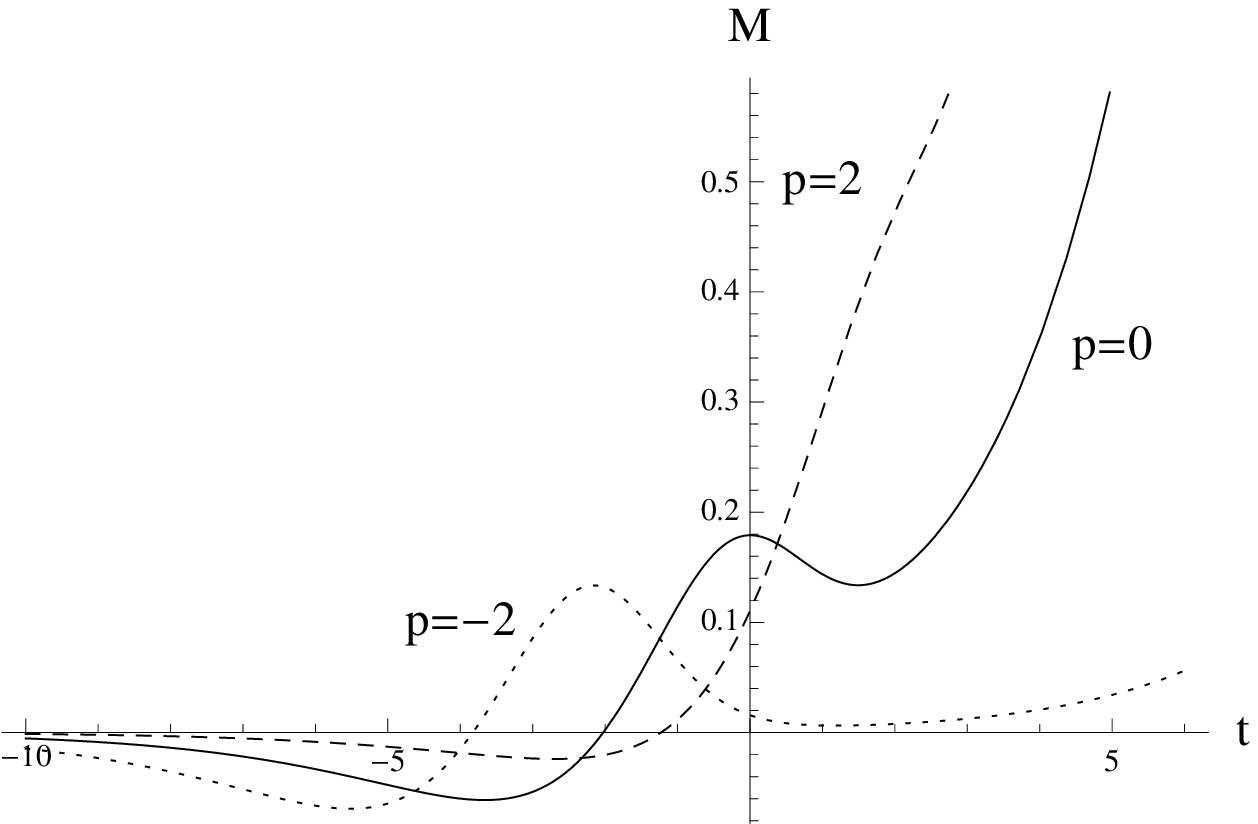}
\includegraphics[width=0.47 \textwidth, height=0.3\textheight]{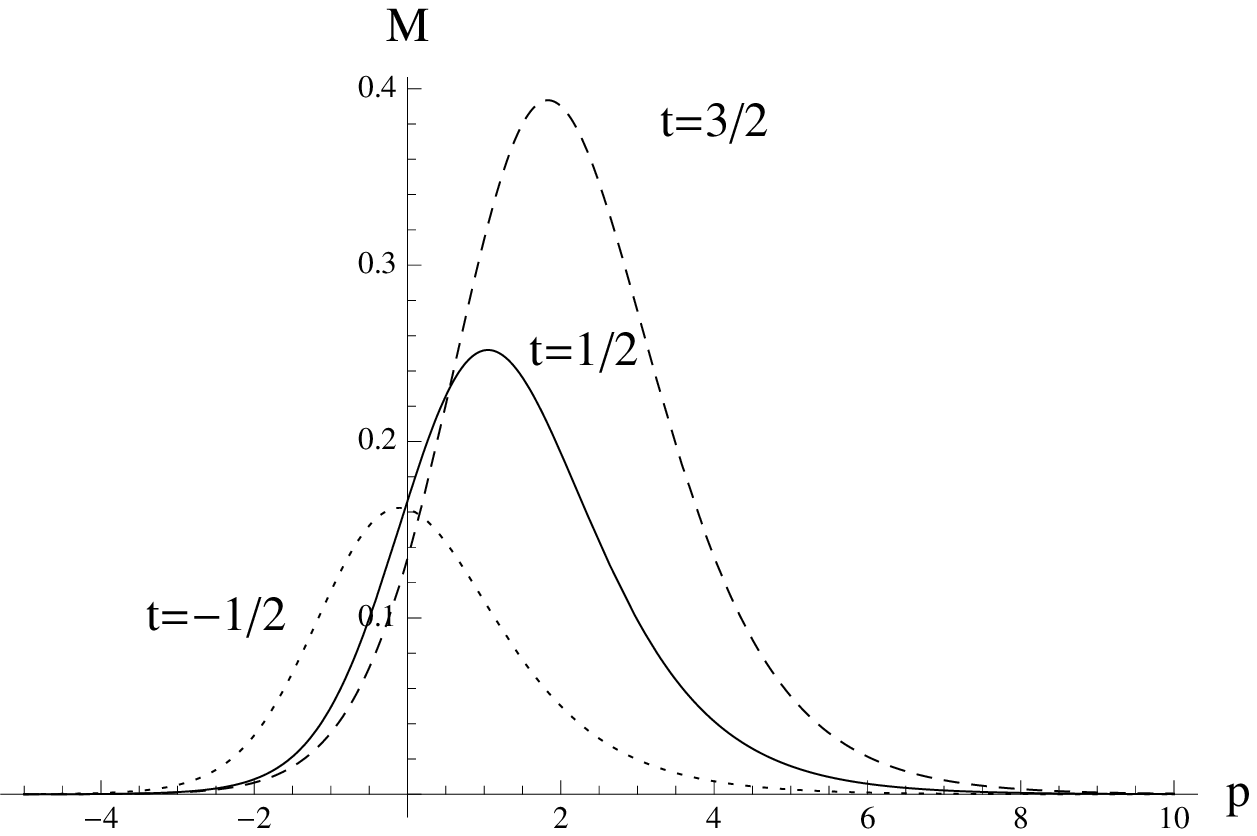}
\hspace*{0.3 \textwidth} (a) \hspace*{0.4 \textwidth} (b) 
\caption{The behaviour of the Melnikov/flux function (\ref{eq:melnikovmadeup}) (a) at different gate positions, and
(b) at different times.}
\label{fig:madeupflux}
\end{figure}

The previous example was area-preserving, enabling the usage of Corollary~\ref{corollary:heteroclinic} in which 
a formal substitution of a Dirac delta impulse into the (smooth) Melnikov function was possible.  For the next example,
consider the non-area-preserving flow given by
\begin{equation}
\fl \frac{d}{d t} \left( \begin{array}{c} x_1 \\ x_2 \end{array} \right) = 
\left( \begin{array}{c} x_1 - x_1^2 \\ 2 x_1 x_2  -\frac{1}{2}  x_2 \end{array} \right) \, ,
\label{eq:madeup}
\end{equation}
whose phase plane is shown in Figure~\ref{fig:madeupphase}.  After an imposed impulsive perturbation, the intention
is to compute the flux across the heteroclinic shown by the thick line, connecting the points $ a \equiv (0,0) $ and 
$ b \equiv (1,0) $.  It is easy to compute that the corresponding heteroclinic trajectory is given by
\[
\fl \bar{x}(t) = \left( \begin{array}{c} \bar{x}_1(t) \\ \bar{x}_2(t) \end{array} \right) = 
\left( \begin{array}{c} \frac{e^{t}}{1 + e^{t}} \\ 0 \end{array} \right) \,  , \, 
f \left( \bar{x}(t) \right) = \left( \begin{array}{c} \frac{e^{t}}{\left(1 + e^{t} \right)^2} \\ 0 \end{array} \right) \, , \, 
f^\perp \left( \bar{x}(t) \right) = \left( \begin{array}{c} 0 \\ \frac{e^{t}}{\left(1 + e^{t} \right)^2} \end{array} \right) \, .
\]
Now if a general impulsive perturbation of the form (\ref{eq:dirac}) is imposed, it has been established that the 
resulting flux, in this case across the heteroclinic from the lower to the upper strips lying within $ 0 < x_1 < 1 $, is given by $ \phi(p,t) = \eps M(p,t) 
+ \O{\eps^2} $ where the Melnikov function is given in (\ref{eq:melnikovdirac}).  In computing this, 
(\ref{eq:jump}) from Theorem~\ref{theorem:distance} yields
\[
\fl j_i(p,t) = \frac{e^{t_i-t+p}}{\left( 1 + e^{t_i-t+p} \right)^2} \, g_{i,2} \left( \frac{ e^{t_i-t+p}}{1 + e^{t_i-t+p}}, 0 \right) \, ,
\]
where $ g_{i,2} $ is the second component of the vector $ g_i $.  In this case, $ {\mathrm{Tr}} \, D f = -2 x_1 + 1 + 2 x_1 -1/2 = 
1/2 $, and the flow is expanding.  Thus $ \hat{F}^{u,s}_p(s) = 1/ (2 s) $, and evaluating both parts of (\ref{eq:resolvent}) 
gives $ R_p(t) = e^{t/2}/2 $ for $ t \ne 0 $.  Therefore, from (\ref{eq:melnikovdirac}), 
\[
\fl M(p,t) = \sum_{i=1}^n \left[ j_i(p,t) + \frac{1}{2} \int_{t_i}^t e^{(t-\xi)/2} j_i(p,\xi) \, \d \xi \right] \, .
\]
Some sample calculations for the specific choice of $ n = 2 $, $ t_1 = 0 $, $ t_2 = 1 $, $ g_{1,2}(x_1,x_2) = e^{t_1} x_1 $
and $ g_{2,2}(x_1,x_2) = x_2^2 + x_1^3 $ are now performed.  Then,
\begin{eqnarray}
\fl M(p,t) & = & \frac{e^{-t+p}}{\left(1 + e^{-t+p} \right)^2} e^0 \frac{e^{-t+p}}{1+e^{-t+p}} + \frac{1}{2} \int_0^t
e^{(t-\xi)/2} \frac{e^{-2 \xi+2p}}{\left(1 + e^{-\xi+p} \right)^3} \, \d \xi \nonumber \\
& & +  \frac{e^{1-t+p}}{\left(1 + e^{1-t+p} \right)^2} \left( \frac{e^{1-t+p}}{1+e^{1-t+p}} \right)^3 + \frac{1}{2} \int_1^t
e^{(t-\xi)/2} \frac{e^{4(1-\xi+p)}}{\left(1+e^{1-\xi+p}\right)^5}  \, \d \xi \, .
\label{eq:melnikovmadeup}
\end{eqnarray}
The above can be explicitly integrated, leading to a not particularly illuminating lengthy expression.  Its behaviour with $ p $ and $ t $ is shown
in Figure~\ref{fig:madeupflux}, bearing in mind that  positive $ M $ relates to flow across the heteroclinic from the lower
to the upper strip.  In (a), $ M $'s $ t $-variation is shown
for several different gate choices $ p $ ($ p = 0 $ would be the midpoint, $ x_1 = 1/2 $).  The flow profile through the
leftmost ($ p = -2 $) gate is seen to gradually flow through the next gates, but with additional accummulating
effects.  As time progresses, it appears that the flux increases without bound, which is unsurprising because the flow is
expanding ($ {\mathrm{Tr}} \, D f > 0 $), and  $ \O{\eps} $-theory is only valid for  $ t \in [T_s,T_u] $.   The rapid decay of
the flux in backwards time is because the flow is compressing in backwards time.  In (b), 
$ M $'s variation with a continuously moving gate $ p $ is shown at several $ t $ values: one below $ t_1 $, one 
between $ t_1 $ and $ t_2 $, and one after $ t_2 $.   The peak flux location 
is initially towards the left, but moves towards the right as time progresses.  

\section{Concluding remarks}
\label{sec:conclude}

This article has formulated and characterised how stable and unstable manifolds in two-dimensional flows are
influenced by state-dependent impulsive perturbations, by casting the problem as an integral equation.  The particular framework
chosen here is associated with thinking of Dirac impulses as a limit of (potentially asymmetric) rectangular pulses.
The methodology
allows for the determination of a condition for persisting heteroclinic connections, and also quantifying fluid flux
across previously impermeable heteroclinic manifolds.  The spatial variation of the impulses, and compressibility of the
flow, are both taken into account.

Extensions to these results are currently being pursued on several fronts.  The numerical difficulties of inverting the
Laplace transform are well-known \cite{kuhlman}; approaches to characterise the pseudo-manifold locations using different formul\ae{} 
would be of value.  Extending the results to higher dimensions, in particular three \cite{3d}, would also be beneficial, since
fluid transport across time-varying two-dimensional surfaces has a profound impact on geophysical and microfluidic
mixing.  The question as to whether, analogous to recent work \cite{saddlecontrol,manifoldcontrol,controlnd}, a Melnikov
approach can be used to control stable and unstable manifolds but in a {\em discontinuous} fashion, is another future
direction of research.  The ability to reformulate the results for impulses which are randomly chosen (e.g., a randomly kicked
Duffing oscillator), to extend to a countable number of impulses, or to formulate the problem for {\em general} 
$ \delta $-families, would also be  of interest.

\ack
The author was supported by the Australian Research Council through Future Fellowship grant FT130100484.  Thanks 
is also expressed to an anonymous referee, whose detailed comments on determining a jump map associated with the
impulsive differential equations approach, led to the discussion in Section~\ref{sec:integral}.

\appendix

\section{Proof of Theorem~\ref{theorem:unstabledirac} (Unstable pseudo-manifold)}
\label{sec:melnikov_impulsive_unstable}

Consider a fixed time-slice $ t $ in the augmented $ \Omega \times (-\infty,T_u] \setminus {\mathcal J} $ phase space,
and a fixed $ p \in (-\infty,P] $.  Let $ \tau $ represent the time-variation henceforth in this proof, since $ t $ is assumed fixed.
Now, when $ \eps = 0 $,
the trajectory $ \bar{x}^u(\tau-t+p) $ is a solution to (\ref{eq:dirac}) such that this passes
through the point $ \bar{x}^u(p) $ in the time-slice $ t $.  When $ \eps \ne 0 $, suppose $ x_\eps^u(p,\tau) $ is a nearby
trajectory lying on the unstable pseudo-manifold, which can be represented by
\begin{equation}
\fl x_\eps^u(p,\tau) = \bar{x}^u(\tau-t+p) + \eps x_1(p,\tau,\eps) \, .
\label{eq:xudirac}
\end{equation}
The quantity $ x_1 $ is $ \O{\eps} $ for $ (p,\tau) \in (-\infty,P] \times (-\infty,T_u] \setminus {\mathcal J} $ since
the effect of $ \eps $ is only to introduce a finite number of $ \O{\eps} $ jumps in the solution.  Moreover, since the
perturbation will only begin affecting solutions for $ t > t_1 $, $ x_1(p,\tau,\eps) $ must be zero for $ \tau < t_1 $.
With the replacement $ \beta \rightarrow - \infty $ and $ t \rightarrow \tau $ in  (\ref{eq:dirac}), the evolution of this
trajectory on the unstable pseudo-manifold satisfies
\begin{equation}
\fl x_\eps^u(p,\tau) = a + \! \int_{-\infty}^\tau  \! \! \! f \left( x_\eps^u(p,\xi) \right)~ \d \xi
+ \eps \sum_{i=1}^n
\I_{(-\infty,\tau)}(t_i)
\left[ \alpha g_i \! \left( x_\eps^u(p,t_i^-) \right) + (1-\alpha) g_i \! \left( x_\eps^u(p,t_i^+) \right) \right] \, .
\label{eq:diracu}
\end{equation}
The $ \eps = 0 $ restriction of (\ref{eq:diracu}), in conjunction with (\ref{eq:xudirac}) indicates that
\begin{equation}
\fl \bar{x}^u(\tau-t+p) = a + \int_{-\infty}^\tau f \left( \bar{x}^u(\xi-t+p)  \right) \, \d \xi \, ,
\label{eq:xuunpdirac}
\end{equation}
which is a simple statement that  $ \bar{x}^u(\tau-t+p) $ obeys the differential equation
\[
\fl \frac{\partial}{\partial \tau} \bar{x}^u(\tau-t+p) = f \left( \bar{x}^u(\tau-t+p) \right)  
\]
and satisfies $  \bar{x}^u(\tau-t+p) \rightarrow a $ as $ \tau \rightarrow - \infty $ for any $ p $ and $ t $.  What is 
required is the $ \O{\eps} $-modification to this solution lying on the unstable manifold when $ \eps \ne 0 $. 
Substituting (\ref{eq:xudirac}) into (\ref{eq:diracu}) gives
\begin{eqnarray*}
\fl \bar{x}^u(\tau-t+p) + \eps x_1(p,\tau,\eps) & = & a + \int_{-\infty}^\tau f \left( \bar{x}^u(\xi-t+p)+\eps x_1(p,\xi,\eps) \right) \, \d \xi \\
& & + \eps \alpha  \sum_{i=1}^n \I_{(-\infty,\tau)}(t_i) g_i \left( \bar{x}^u(t_i^- -t+p) + \eps x_1(p,t_i^-,\eps) \right) \\
& & + \eps (1-\alpha) \sum_{i=1}^n \I_{(-\infty,\tau)}(t_i) g_i \left( \bar{x}^u(t_i^+ -t+p) + \eps x_1(p,t_i^+,\eps)  \right) \, .
\end{eqnarray*}
Next, $ f $ and each $ g_i $ will be Taylor expanded around $ \bar{x}^u $.  Terms beyond
$ \O{\eps} $ will include $ D^2 f $ and $ D g_i $, all of which are bounded on $ \Omega \times \R $ by Hypotheses~\ref{hyp:fdirac}
and \ref{hyp:gdirac}.  While the $ D g_i $ terms appear in a regular fashion, the $ D^2 f $ terms appear inside an integral 
over an unbounded domain, but  since the $ x_1(p,\tau,\eps) $ appearing in the integrand is zero
for $ \tau < t_1 $, all these terms are $ \O{\eps^2} $. Thus,
\begin{eqnarray*}
\fl \bar{x}^u(\tau-t+p) + \eps x_1(p,\tau,\eps) & = & a + \int_{-\infty}^\tau f \left( \bar{x}^u(\xi-t+p) \right) \, \d \xi \\
& & + \int_{-\infty}^\tau D f \left( \bar{x}^u(\xi-t+p) \right) \eps x_1(p,\xi,\eps) \, \d \xi + \O{\eps^2} \\
& & \hspace*{-2cm} +  \eps \sum_{i=1}^n \I_{(-\infty,\tau)}(t_i) \left[ \alpha g_i \left( \bar{x}^u(t_i^- -t+p) \right) +
(1-\alpha) g_i \left( \bar{x}^u(t_i^+ -t+p) \right) \right] \, .
\end{eqnarray*}
Utilising (\ref{eq:xuunpdirac}) and the continuity of $ \bar{x}^u $ and $ g_i $,
\begin{equation}
\fl x_1(p,\tau,\eps) = \! \! \int_{-\infty}^\tau \! \! \! D f \left( \bar{x}^u(\xi \! - \! t \! + \! p) \right) x_1(p,\xi,\eps) \, \d \xi +  \sum_{i=1}^n \I_{(t_i,\infty)} \! (\tau)
g_i \! \left( \bar{x}^u\left(t_i \! - \! t \! + \! p \right) \right) + \O{\eps} \, .
\label{eq:x1dirac}
\end{equation}
Of interest is the fact that in this leading-order expression for the unstable pseudo-manifold, the $ \alpha $-dependence
has dropped out.  Impacts on asymmetry of the Dirac impulse representation are therefore only felt at higher-order.
Now define
\begin{eqnarray}
\fl \tilde{M}^u(p,\tau,\eps) & := & f \left( \bar{x}^u(\tau-t+p) \right)^\perp \cdot \frac{ x_\eps^u(p,\tau) - \bar{x}^u(\tau-t+p)}{\eps} \nonumber \\
& = & 
f \left( \bar{x}^u(\tau-t+p) \right)^\perp \cdot x_1(p,\tau,\eps) \, , 
\label{eq:mudiracdef}
\end{eqnarray}
and notice that $ \tilde{M}^u(p,t,\eps) = f \left( \bar{x}^u(p) \right)^\perp x_1(p,t,\eps) $ 
expresses the leading-order displacement of the unstable pseudo-manifold in the normal direction to the original manifold
at a point $ \bar{x}^u(p) $, in the time-slice $ t $.  
Given that $ x_1(p,-\infty,\eps) = 0 $ and $ \tilde{M}^u(p,-\infty,\eps) = 0 $, it is possible to rewrite (\ref{eq:mudiracdef}) in the form
\[
\fl \tilde{M}^u(p,\tau,\eps) =  \int_{-\infty}^\tau \frac{d}{d \xi} \left[ f \left( \bar{x}^u(\xi-t+p) \right)^\perp \cdot x_1(p,\xi,\eps) \right] \,
\d \xi \, . 
\]
In writing the above, it has been noted that while $ x_1(p,\xi,\eps) $ is not differentiable in $ \xi $ at the jump values $ t_i $, its
temporal derivative {\em is} integrable.  By using the product rule in the integrand,
\begin{eqnarray}
\fl \tilde{M}^u(p,\tau,\eps) & = & \int_{-\infty}^\tau \frac{d}{d \xi} \left[ f \left( \bar{x}^u(\xi-t+p) \right)^\perp \right] \cdot 
 x_1(p,\xi,\eps) \, \d \xi \nonumber \\
 & & + \int_{-\infty}^\tau f \left( \bar{x}^u(\xi-t+p) \right)^\perp \cdot 
 \frac{d}{d \xi} \left[  x_1(p,\xi,\eps) \right] \, \d \xi \nonumber\\
 & = & \int_{-\infty}^\tau \left[ D f \left( \bar{x}^u(\xi-t+p) \right) f \left( \bar{x}^u(\xi-t+p) \right) \right]^\perp \cdot \! 
 x_1(p,\xi,\eps) \, \d \xi + \O{\eps} \nonumber \\
 & & \hspace*{-3cm} + \int_{-\infty}^\tau \! \! \! f \left( \bar{x}^u(\xi \! - \! t \! + \! p) \right)^\perp \! \! \cdot \! \left[ \! 
 D f \left( \bar{x}^u(\xi \! -\! t \! + \! p) \right) x_1(p,\xi,\eps) \! + \!  
 \sum_{i=1}^n  g_i \left( \bar{x}^u(t_i \! -\! t \! + \! p) \right)  \! \frac{d}{d\xi}\I_{(t_i,\infty)}(\xi) \right] \d \xi \nonumber \\
  \nonumber
\end{eqnarray}
where the second equality is by taking the derivative of (\ref{eq:x1dirac}) in a distributional sense.  Now, using
easily verifiable identity for vectors $ b $ and $ c $ in $ \R^2 $, and $ 2 \times 2 $ matrices $ A $, given by (see \cite[e.g.]{tangential,holmes})
\[
\left( A b \right)^\perp \cdot c + b^\perp \cdot \left( A c \right) = {\mathrm{Tr}} A \left( b^\perp \cdot c \right) \, ,
\]
and by choosing $ A = D f $, $ b = f $ and $ c = x_1 $, 
\begin{eqnarray}
 \fl \tilde{M}^u(p,\tau,\eps) & = & \int_{-\infty}^\tau {\mathrm{Tr}} \, D f \left( \bar{x}^u(\xi-t+p) \right) f \left( \bar{x}^u(\xi-t+p) \right)^\perp \cdot x_1(p,\xi,\eps)
 \, \d \xi + \O{\eps} 
\nonumber \\
 & & \hspace*{-2cm} + \sum_{i=1}^n  g_i \left( \bar{x}^u(t_i-t+p)  \right) \cdot 
 \int_{-\infty}^\tau f \left( \bar{x}^u(\xi-t+p) \right)^\perp \frac{d}{d\xi}\I_{(t_i,\infty)}(\xi)  \d \xi \, .
 \label{eq:mudiractemp}
 \end{eqnarray}
The distributional integral above is now evaluated by parts:
\begin{eqnarray*}
\fl  \int_{-\infty}^\tau f \left( \bar{x}^u(\xi-t+p) \right)^\perp \frac{d}{d\xi}\I_{(t_i,\infty)}(\xi)  \d \xi 
& = &  f \left( \bar{x}^u(\xi-t+p) \right)^\perp \I_{(t_i,\infty)}(\xi) \Big]_{\xi=-\infty}^\tau \\
& & - \int_{-\infty}^\tau \frac{d}{d\xi} \left[
 f \left( \bar{x}^u(\xi-t+p) \right)^\perp \right]  \I_{(t_i,\infty)}(\xi) \d \xi \\
 & = &  \left[ f \left( \bar{x}^u(\tau-t+p) \right)^\perp \I_{(t_i,\infty)}(\tau) - 0 \right] \\
 & & - \I_{(t_i,\infty)}(\tau) \int_{t_i}^\tau
 \frac{d}{d\xi} \left[  f \left( \bar{x}^u(\xi-t+p) \right)^\perp \right] \d \xi \\
 & = &  f \left( \bar{x}^u(\tau-t+p) \right)^\perp \I_{(t_i,\infty)}(\tau) \\
 & & - \I_{(t_i,\infty)}(\tau) \! \left[ f \left( \bar{x}^u(\tau \! - \! t \! + \! p) \right)^\perp \! \! 
 - \! f \left( \bar{x}^u(t_i \! - \! t \! + \! p) \right)^\perp  \right] \\
 & = & \I_{(t_i,\infty)}(\tau)  f \left( \bar{x}^u(t_i-t+p) \right)^\perp \, .
 \end{eqnarray*}
Using this and the fact that $ \tilde{M}^u = f^\perp \cdot x_1 $ into (\ref{eq:mudiractemp}) gives 
\begin{eqnarray}
\fl \tilde{M}^u(p,\tau,\eps) & =  & \int_{-\infty}^\tau {\mathrm{Tr}} \, D f \left( \bar{x}^u(\xi-t+p) \right) \tilde{M}^u(p,\xi,\eps)
 \, \d \xi + \O{\eps} \nonumber \\
 & & + \sum_{i=1}^n \I_{(t_i,\infty)}(\tau)  f \left( \bar{x}^u(t_i-t+p) \right)^\perp \cdot g_i \left( \bar{x}^u(t_i-t+p) \right) \, .
 \label{eq:mutildetemp}
 \end{eqnarray}
 Let $ M^u(p,t) $ be the solution to the above with the $ \O{\eps} $ term neglected, which also satisfies $ M^u(p,\tau) = 0 
 $ for $ \tau < t_1 $.  Then.
 \[
\fl  \left[ \tilde{M}^u(p,\tau,\eps) - M^u(p,\tau) \right]  = \int_{t_1}^\tau {\mathrm{Tr}} \, D f \left( \bar{x}^u(\xi-t+p) \right) 
  \left[ \tilde{M}^u(p,\tau,\eps) - M^u(p,\tau) \right] \, \d \xi + \O{\eps} \, , 
  \]
 and thus $ \tilde{M}^u(p,\tau,\eps) - M^u(p,\tau) = \O{\eps} $ for $ \tau \in (-\infty,T_u] \setminus {\mathcal J} $.
Therefore, the interchange of  $ \tilde{M} $ with $ M $ in  (\ref{eq:mudiracdef}) is legitimate, 
since this will only cause a $ \O{\eps^2} $ error in the normal distance
measure.  So replacing $ \tilde{M} $ with $ M $, neglecting the $ \O{\eps} $ term, and replacing $ \tau $ above with $ t $ 
 leads to the integral equation for the unstable
Melnikov function:
\begin{eqnarray}
\fl M^u(p,t) & = & \int_{-\infty}^t {\mathrm{Tr}} \, D f \left( \bar{x}^u(\xi-t+p) \right) M^u(p,\xi)
 \, \d \xi  \nonumber \\
 & & + \sum_{i=1}^n \I_{(t_i,\infty)}(t)  f \left( \bar{x}^u(t_i-t+p) \right)^\perp \cdot g_i \left( \bar{x}^u(t_i-t+p) \right) \nonumber \\
 & = & \int_{-\infty}^t {\mathrm{Tr}} \, D f \left( \bar{x}^u(\xi \! - \! t \! + \! p) \right) M^u(p,\xi) \, \d \xi 
 +  \sum_{i=1}^n \I_{(t_i,\infty)}(t) j_i^u(p,t) \, , 
 \label{eq:muintegraleq}
 \end{eqnarray}
 where the definition (\ref{eq:jumpu}) has been used.  The equation (\ref{eq:muintegraleq}) is a Volterra equation of the
 second kind, but over an unbounded domain, and with a discontinuous inhomogeneity.
The following lemma will help in solving (\ref{eq:muintegraleq}).

\begin{lemma}
\label{lemma:renewal}
Consider the integral equation
\begin{equation}
\fl M(t) = j(t) + \int_{-\infty}^t M(\xi) F(t-\xi) \, d \xi
\label{eq:renewal}
\end{equation}
where $ j $ is piecewise differentiable and is zero below some finite value $ t_1 $, and the kernel $ F \in {\mathrm{C}}^1 \left(
[0,\infty) \right) $ satisfies $ \left| \lim_{t \rightarrow \infty} F(t) \right| = F_0 < \infty $.  Then, (\ref{eq:renewal}) has a solution
\begin{equation}
\fl M(t) = j(t) + \int_{-\infty}^t R(t-\xi) j(\xi) \, \d \xi
\label{eq:renewalsol}
\end{equation}
at values $ t $ at which $ j $ is defined, 
where the resolvent $ R $ is obtained from the Laplace transform $ \hat{F}(s) $ of $ F(t) $ by $ R(t) = {\mathcal L}^{-1} \left\{
\hat{F}(s) / \left[ 1 - \hat{F}(s)\right] \right\}(t) $.
\end{lemma}

\proof
Equation~(\ref{eq:renewal}) is in the form of a renewal equation \cite{polyaninmanzhirov,bellmancooke} but with an unbounded domain.  The basic
renewal equation solution with the infinite limit substituted is indeed (\ref{eq:renewalsol}) \cite{polyaninmanzhirov,bellmancooke}.
However, the legitimacy of this formal process requires the conditions on  $ F $ and $ j $ as given in Lemma~\ref{lemma:renewal},
based on which a full proof is given in \ref{sec:melnikov_impulsive_renewal}.
\proofend

To now prove Theorem~\ref{theorem:unstabledirac}, Lemma~\ref{lemma:renewal} is applied to (\ref{eq:muintegraleq}) with the choice
\[
\fl F(t) = {\mathrm{Tr}} \, D f \left( \bar{x}^u(p-t) \right) \quad {\mathrm{and}} \quad
j(t) = \sum_{i=1}^n \I_{(t_i,\infty)}(t) j_i^u(p,t) \, .
\]
Using the resolvent definition (\ref{eq:resolventu}), this yields
\[
\fl M^u(p,t) = \sum_{i=1}^n \I_{(t_i,\infty)}(t) j_i^u(p,t)
 + \int_{-\infty}^t R_p^u(t-\xi)  \sum_{i=1}^n \I_{(t_i,\infty)}(\xi) j_i^u(p,\xi) \, \d \xi \, , 
\]
from which (\ref{eq:mudirac}) arises since each of the jump functions is only turned on for $ \xi $ values greater than
$ t_i $ in the integrand. Thereby, Theorem~\ref{theorem:unstabledirac} has been proven.

\subsection{Proof of Lemma~\ref{lemma:renewal} (Integral equation for $ M^u $)}
\label{sec:melnikov_impulsive_renewal}

In proving Lemma~\ref{lemma:renewal}, a preliminary lemma proves convenient.

\begin{lemma}
\label{lemma:renewaltemp}
Let $ F $ and $ j $ satisfy the hypotheses stated in Lemma~\ref{lemma:renewal}.  If $ w(t) $ satisfies
\begin{equation}
w(t) = 1 + \int_0^t w(\xi) F(t-\xi) \, \d \xi \, , 
\label{eq:renewaltemp}
\end{equation}
for $ t \ge 0 $, then the solution to the integral equation (\ref{eq:renewal}) is given by 
\begin{equation}
M(t) = \int_{-\infty}^t w(t-\xi) \frac{d}{d \xi} \left[ j(\xi) \right] \, \d \xi \, , 
\label{eq:renewalsoltemp}
\end{equation}
where since $ j $ is piecewise continuous the derivative in (\ref{eq:renewalsoltemp}) is to be considered in a 
distributional sense.
\end{lemma}

\proof
Define the potential solution
\[
\fl \bar{M}(t) :=  \int_{-\infty}^t w(t-\xi) \frac{d}{d \xi} \left[ j(\xi) \right] \, \d \xi \, .
\]
Now, closely following \cite{bellmancooke}, 
\begin{eqnarray*}
\fl \int_{-\infty}^t \bar{M}(\eta) F(t-\eta) \d \eta & = & \int_{-\infty}^t \left[ \int_{-\infty}^\eta w(\eta-\xi) \frac{d}{d \xi} \left[ j(\xi) \right]
\d \xi \right] F(t-\eta) \, \d \eta \\
& = & \int_{-\infty}^t \frac{d}{d \xi} \left[ j(\xi) \right] \int_{\xi}^t w(\eta-\xi) F(t-\eta) \d \eta \, \d \xi \\
& = & \int_{-\infty}^t  \frac{d}{d \xi} \left[ j(\xi) \right] \int_{0}^{t-\xi} w(u) F \left( \left[ t - \xi \right] - u \right) \, \d u \, \d \xi \\
& = & \int_{-\infty}^t  \frac{d}{d \xi} \left[ j(\xi) \right]  \left[ w(t-\xi) - 1 \right] \, \d \xi  \quad \left[{\mathrm{by~}} (\ref{eq:renewaltemp})
\right] \\
& = & \int_{-\infty}^t w(t-\xi) \frac{d}{d \xi} \left[ j(\xi) \right] \, \d \xi - \int_{-\infty}^t \frac{d}{d \xi} \left[ j(\xi) \right] \d \xi \\
& = & \bar{M}(t) - j(\xi) \Big]_{\xi=-\infty}^t = \bar{M}(t) - j(t) \, ,
\end{eqnarray*}
which proves that  $ \bar{M}(t) $ does indeed satisfy (\ref{eq:renewal}).  The interchanging of the order of integration in 
was legitimate  since the integrand was therefore absolutely integrable over the unbounded domain ($ j(t) = 0 $ for $ t <t_1 $ while $ F(t-\eta) $ approached a limit as $ \eta 
\rightarrow - \infty $).  
\proofend

Lemma \ref{lemma:renewaltemp} reduces the problem to finding a the solution to the auxilliary equation (\ref{eq:renewaltemp}).  Since the functions here are smooth, it is an easier problem.  Taking the Laplace transform
of (\ref{eq:renewaltemp}), along with the identifications 
$ \hat{w}(s) := {\mathcal L} \left\{ w(t) \right\}(s) $ and $ \hat{F}(s):= {\mathcal L} \left\{
F(t) \right\}(s) $ gives
\[
\fl \hat{w}(s) = \frac{1}{s} + \hat{w}(s) \hat{F}(s) \, .
\]
Solving for $ \hat{w}(s) $ gives the result
\[
\fl \hat{w}(s) = \frac{1}{s \left[ 1 - \hat{F}(s) \right]} = \frac{1}{s} + \frac{1}{s} \frac{ \hat{F}(s) }{1 - \hat{F}(s)} =: \frac{1}{s} +
\frac{1}{s} \, \hat{R}(s) \, .
\]
Inverting the Laplace transform and once again using the convolution property gives
\[
\fl w(t) = 1 + \int_0^t R(\xi) \, \d \xi \, .
\]
Inserting the above into (\ref{eq:renewalsoltemp}) then results in
\begin{eqnarray*}
\fl M(t) & = & \int_{-\infty}^t \left( 1 + \int_0^{t-\xi} R(\eta) \, \d \eta \right) \frac{d}{d \xi} \left[ j(\xi) \right] \, \d \xi \\
& = & \int_{-\infty}^t \frac{d}{d \xi} \left[ j(\xi) \right] \, \d \xi + \int_{-\infty}^t \left( \int_0^{t-\xi} R(\eta) \, \d \eta \right) 
\frac{d}{d \xi} \left[ j(\xi) \right] \, \d \xi \\
& = & j(t) + \left[ \left( \int_0^{t-\xi} R(\eta) \d \eta \right) j(\xi) \Big]_{\xi=-\infty}^t - \int_{-\infty}^t \left[ - R(t-\xi) \right]
j(\xi) \, \d \xi \right] \\
& = & j(t) + \int_{-\infty}^t R(t-\xi) j(\xi) \, \d \xi  \, ;
\end{eqnarray*}
the result required for Lemma~\ref{lemma:renewal}.

\section{Proof of Theorem~\ref{theorem:stabledirac} (Stable pseudo-manifold)}
\label{sec:melnikov_impusive_stable}

Details which are similar to, and with obvious modifications from, the proof of the unstable pseudo-manifold expressions
of Theorem~\ref{theorem:unstabledirac} as given
in \ref{sec:melnikov_impulsive_unstable} will be sketched briefly.  
However, there are some issues---in particular dealing with how the Laplace transform representation is to be modified
for functions with negative argument---for which more details will be given.  

As in \ref{sec:melnikov_impulsive_unstable}, consider a fixed time-slice $ t $ and a fixed $ p $, and let $ \tau $
be the time-variable.  Define
\begin{equation}
\fl x_\eps^s(p,\tau) := \bar{x}^s(\tau-t+p) + \eps x_1(p,\tau,\eps)
\label{eq:xsdirac}
\end{equation}
where now $ x_1 $ is $ \O{\eps} $ for $ (p,\tau) \in [P,\infty) \times [T_s,\infty) \setminus {\mathcal J} $.  Define also
\begin{equation}
\fl M^s(p,\tau,\eps):= 
f \left( \bar{x}^s(\tau-t+p) \right)^\perp \cdot x_1(p,\tau,\eps) \, . 
\label{eq:msdiracdef}
\end{equation}
Now, consider using
the evolution equation (\ref{eq:dirac}) with $ \beta = \infty $.  Following a derivation similar to \ref{sec:melnikov_impulsive_unstable},  instead of (\ref{eq:muintegraleq}) the integral equation
\begin{equation}
\fl M^s(p,t) =  - \int_t^\infty {\mathrm{Tr}} \, D f \left( \bar{x}^s(\xi-t+p) \right) M^s(p,\xi)
 \, \d \xi  - \sum_{i=1}^n \I_{(-\infty, t_i)}(t) j_i^s(p,t)
 \label{eq:msintegraleq}
 \end{equation}
results for the stable Melnikov function $ M^s $, where $ j_i^s $ is defined in (\ref{eq:jumps}).   This integral equation can be 
solved with the help of the following lemma, analogous to Lemma~\ref{lemma:renewal}.

\begin{lemma}
\label{lemma:renewals}
Consider the integral equation
\begin{equation}
\fl M(t) = -j(t) - \int_t^\infty M(\xi) F(t-\xi) \, d \xi
\label{eq:renewals}
\end{equation}
where $ j $ is piecewise differentiable and is zero above some finite value $ t_n $, and the kernel $ F \in {\mathrm{C}}^1 \left( 
(-\infty,0] \right) $ satisfies $ \left| \lim_{t \rightarrow - \infty} F(t) \right| = F_0 < \infty $.  Then, (\ref{eq:renewals}) has a solution
\begin{equation}
\fl M(t) = -j(t) + \int_t^\infty R(t-\xi) j(\xi) \, \d \xi
\label{eq:renewalsols}
\end{equation}
at values $ t $ at which $ j $ is defined, 
where the resolvent $ R $ is obtained from the Laplace transform $ \hat{F}(s) $ of $ F(-t) $ by $ R(t) = {\mathcal L}^{-1} \left\{
\hat{F}(s) / \left[ 1 + \hat{F}(s)\right] \right\}(-t) $.
\end{lemma}

\proof
The first claim is that if $ w(t) $ solves
\begin{equation}
\fl w(t) = -1 + \int_0^t w(\xi) F(t-\xi) \, \d \xi \, , 
\label{eq:renewaltemps}
\end{equation}
for $ t \le 0 $ (with $ w $ and $ F $ being defined on $ (0,\infty) $), then the solution to the integral equation (\ref{eq:renewals}) is given by 
\begin{equation}
\fl M(t) = - \int_t^\infty w(t-\xi) \frac{d}{d \xi} \left[ j(\xi) \right] \, \d \xi \, . 
\label{eq:renewalsoltemps}
\end{equation}
The proof of this is similar to that of Lemma~\ref{lemma:renewaltemp} and will be skipped.  To use Laplace transform methods
to solve (\ref{eq:renewaltemps}), replacing $ t $ with $ - t $ enables the representation
\[
\fl w(-t) = - 1 + \int_0^{-t} w(\xi) F(-t-\xi) \, \d \xi 
\]
with domain of validity now $ t \ge 0 $.
Defining $ \tilde{w}(t) = w(-t) $ and $ \tilde{F}(t) = F(-t) $ results in 
\[
\fl \tilde{w}(t) = -1 + \int_0^{-t} \tilde{w}(-\xi) \tilde{F}(t + \xi) \, \d \xi = -1 - \int_0^{t} \tilde{w}(\eta) \tilde{F}(t -\eta) \, \d \eta \, .
\]
Since each of $ \tilde{w} $ and $ \tilde{F} $ are defined for
$ t \ge 0 $, it is possible to define $ \hat{F}(s) = {\mathcal L} \left\{ \tilde{F}(t) \right\}(s) = 
{\mathcal L} \left\{ F(-t) \right\}(s) $ and $ \hat{w}(s) = {\mathcal L} \left\{ \tilde{w}(t) \right\}(s) = {\mathcal L} \left\{
w(-t) \right\}(s) $.  Taking the Laplace transform of the above expression gives
\[
\fl \hat{w}(s) = - \frac{1}{s} - \hat{w}(s) \hat{F}(s) \, , 
\]
and therefore
\[
\fl \hat{w}(s) = - \frac{1}{s} + \frac{1}{s} \frac{\hat{F}(s)}{1+\hat{F}(s)} \, .
\]
Let $ \hat{R}(s) = \hat{F}(s) / \left[ 1 + \hat{F}(s) \right] $, with inverse Laplace transform $ \tilde{R}(t) $, which is defined for
$ t \ge 0 $.  Applying the
convolution property  yields
\[
\fl \tilde{w}(t) = -1 + \int_0^t \tilde{R}(\xi) \, \d \xi
\]
which with the replacement $ t \rightarrow - t $ gives
\[
\fl w(t) = -1 + \int_0^{-t} \tilde{R}(\xi) \, \d \xi
\]
where now $ t \le 0 $.  This solution for $ w $ when inserted into (\ref{eq:renewalsoltemps}) yields
\begin{eqnarray*}
\fl M(t) & = & - \int_t^\infty \left[ -1 + \int_0^{-t+\xi} \tilde{R}(\eta) \, \d \eta \right] \frac{d}{d \xi} \left[ j(\xi) \right] \, \d \xi \\
& = & -j(t) - \left[  \left( \int_0^{-t+\xi} \tilde{R}(\eta) \, d \eta \right) j(\xi) \Big]_{\xi = t}^\infty - \int_t^\infty \tilde{R}(-t+\xi) j(\xi) \, \d \xi \right] \\
& = & -j(t) + \int_t^\infty \tilde{R}(-t+\xi) j(\xi) \, \d \xi \\
& = & -j(t) + \int_t^\infty R(t-\xi) j(\xi) \, \d \xi 
\end{eqnarray*}
where $ R(t) := \tilde{R}(-t) $ was used to express the solution in terms of a resolvent $ R $ defined for $ t \le 0 $.  
This is the result required.
\proofend

The result of Lemma~\ref{lemma:renewals} can now be applied to the integral equation (\ref{eq:msintegraleq}) with the choice $ M(t) = M^s(p,t) $, 
$ F(t) = {\mathrm{Tr}} \, D f \left( \bar{x}^s(p-t) \right) $ (and hence $ \tilde{F}(t) = {\mathrm{Tr}} \, D f \left( \bar{x}^s(p+t) \right) $,
whose Laplace transform is defined for $ t \ge 0 $)
and $ j(t) = \I_{(-\infty,t_i)}(t) j_p^s(t) $, to yield
\[
\fl M^s(p,t) = - \sum_{i=1}^n \I_{(-\infty,t_i)}(t) j_i^s(p,t)+ \int_t^\infty R_p^s(t-\xi)  \sum_{i=1}^n \I_{(-\infty,t_i)}(\xi) j_i^s(p,\xi) \, \d \xi \, .
\]
Restricting the integral in relation to the indicator functions gives 
the result of Theorem~\ref{theorem:stabledirac}.

\section{Proof of Theorem~\ref{theorem:distance} (Distance between pseudo-manifolds)}
\label{sec:distance}

Consider the point $ x_\eps^u(p,t) $ which lies on $ \tilde{\Gamma}_\eps^u(a) $ but is along the normal vector to
$ \bar{x}(p) $, as shown in Figure~\ref{fig:distance}.  
From Theorem~\ref{theorem:unstabledirac}, its displacement from $ \bar{x}(p) $
along the normal direction $ \hat{f}^\perp \left( \bar{x}(p) \right) $ is given by $ \eps M^u(p,t) / \left| f \left( \bar{x}(p) \right)
\right| + \O{\eps^2} $, where
\[
\fl M^u(p,t) =  \sum_{i=1}^n \I_{(t_i,\infty)}(t) j_i(p,t) + \sum_{i=1}^{{\mathrm{max}} \left\{ j: t_j < t \right\}} \int_{t_i}^t R_p^u(t-\xi) j_i(p,\xi) \, \d \xi
\]
in which
\[
\fl j_i(p,t) = f^\perp \left( \bar{x}(t_i-t+p) \right) \cdot g_i \left(  \bar{x}(t_i-t+p) \right) \, , 
\]
and $ R_p^u: \R^+ \rightarrow \R $ is defined by
\[
\fl R_p^u(t) = {\mathcal L}^{-1} \left\{ \frac{\hat{F}_p^u(s)}{1 - \hat{F}_p^u(s)} \right\}(t) \quad , \quad
\hat{F}_p^u(s) := {\mathcal L} \left\{ {\mathrm{Tr}} \, D f \left( \bar{x}(p-t) \right) \right\}(s) \, .
\]
Similarly from Theorem~\ref{theorem:stabledirac}, the point $ x_\eps^s(p,t) $ in Figure~\ref{fig:distance} has
a displacement from $ \bar{x}(p) $ in the normal direction given by $ \eps M^s(p,t) / \left| f \left( \bar{x}(p) \right)
\right| + \O{\eps^2} $, in which
\[
\fl M^s(p,t) =  - \sum_{i=1}^n \I_{(-\infty,t_i)}(t) j_i(p,t) + \sum_{i={\mathrm{min}} \left\{ j: t_j > t \right\}}^n \int_t^{t_i} R_p^s(t-\xi) j_i(p,\xi) \, \d \xi
\]
where $ R_p^s: \R^- \rightarrow \R $ is defined by
\[
\fl R_p^u(t) = {\mathcal L}^{-1} \left\{ \frac{\hat{F}_p^s(s)}{1 + \hat{F}_p^s(s)} \right\}(-t) \quad , \quad
\hat{F}_p^s(s) := {\mathcal L} \left\{ {\mathrm{Tr}} \, D f \left( \bar{x}(p+t) \right) \right\}(s) \, .
\]
Letting $ \hat{N}(p) = \hat{N}^{s,u}(p) $ (since they are identical), 
\[
\fl \left[ x_\eps^u(p,t) - x_\eps^s(p,t) \right] \cdot \hat{N}(p)  =  \eps \frac{M^u(p,t) - M^s(p,t)}{\left|
f \left( \bar{x}(p) \right) \right|} + \O{\eps^2} =:  \eps \frac{M(p,t)}{\left|
f \left( \bar{x}(p) \right) \right|} + \O{\eps^2}
\]
with the definition $ M(p,t) = M^u(p,t) - M^s(p,t) $, which can be further simplied according to
\begin{eqnarray*}
\fl M(p,t) & = & M^u(p,t) - M^s(p,t) \\
& = &  \sum_{i=1}^n \left[ \I_{(t_i,\infty)}(t) j_i(p,t)  + \I_{(-\infty,t_i)}(t) j_i(p,t) \right] \\
& & + \sum_{i=1}^{{\mathrm{max}} \left\{ j: t_j < t \right\}} \int_{t_i}^t R_p^u(t-\xi) j_i(p,\xi) \, \d \xi 
-  \sum_{i={\mathrm{min}} \left\{ j: t_j > t \right\}}^n \int_t^{t_i} R_p^s(t-\xi) j_i(p,\xi) \, \d \xi \\
& = &  \sum_{i=1}^n  j_i(p,t) + \sum_{i=1}^{{\mathrm{max}} \left\{ j: t_j < t \right\}} \int_{t_i}^t R_p(t-\xi) j_i(p,\xi) \, \d \xi \\
& & +  \sum_{i={\mathrm{min}} \left\{ j: t_j > t \right\}}^n \int_{t_i}^t R_p(t-\xi) j_i(p,\xi) \, \d \xi \\
& = &  \sum_{i=1}^n  j_i(p,t) + \sum_{i=1}^n \int_{t_i}^t R_p(t-\xi) j_i(p,\xi) \, \d \xi \, , 
\end{eqnarray*}
since $ R_p(t-\xi) $ is equal to $ R_p^u(t-\xi) $ in the first integrand since $ t > \xi $, and to $ R_p^s(t-\xi) $ in 
the second since 
$ t < \xi $.

\mbox{} \hrule \mbox{}

%

\end{document}